\documentclass[12pt]{article}
\usepackage{amssymb}
\usepackage{amsmath}
\usepackage{amsthm}
\usepackage{color}
\usepackage{comment}
\usepackage{geometry}		
\usepackage{graphicx}

\DeclareFontFamily{OT1}{pzc}{}
\DeclareFontShape{OT1}{pzc}{m}{it}{<-> s * [1.10] pzcmi7t}{}
\DeclareMathAlphabet{\mathpzc}{OT1}{pzc}{m}{it}

\let\originalleft\left
\let\originalright\right
\renewcommand{\left}{\mathopen{}\mathclose\bgroup\originalleft}
\renewcommand{\right}{\aftergroup\egroup\originalright}

\newtheorem{theorem}{Theorem}[section]

\newtheorem{corollary}[theorem]{Corollary}
\newtheorem{lemma}[theorem]{Lemma}
\newtheorem{proposition}[theorem]{Proposition}

\theoremstyle{definition}
\newtheorem{definition}{Definition}[section]

\theoremstyle{remark}

\begin{document}

\title{The boundary of chaos for piecewise smooth maps and the boundary of positive Hausdorff dimension for survivor sets of open maps}

\author{Paul Glendinning\begin{footnote}{email: p.a.glendinning@manchester.ac.uk}\end{footnote} \ and Cl\'ement Hege\\
Department of Mathematics, University of Manchester, \\
Manchester M13 9PL, U.K.}
\maketitle

\begin{abstract}
We describe the boundary of chaos separating regions of parameter space with positive topological entropy from those with zero
topological entropy for a class of piecewise smooth maps. This coincides with the boundary of positive Hausdorff dimension for the survivor sets of a class of open maps. There are precisely two types of codimension one transitions across the boundaries. One of these involves heteroclinic connections, and in this case there is a finite number of periodic orbits at the transition point. The other involves an infinite cascade of bifurcations creating infinitely many periodic orbits on the boundary in a sequence called the anharmonic cascade. 
\end{abstract}


\vspace{2pc}

\maketitle


\section{Introduction}\label{sect:intro}Understanding transitions to chaos, i.e. how systems evolve from non-chaotic behaviour (having zero topological entropy) to chaos (having positive topological entropy) as parameters vary, is one of the key problems in 
bifurcation theory. The pioneering work of Feigenbaum \cite{Feigenbaum1978} and Tresser and Coullet \cite{Tresser1978} 
established the period-doubling route to chaos as a robust transition for one dimensional maps of the interval; the robustness property is emphasised in \cite{JSY2010}. For maps of the circle there are further transitions including circle intermittency and a heteroclinic transition \cite{MT1986}. The existence of discontinuities can complicate the picture. Glendinning \cite{G1992,G1993} describes a transition for piecewise monotonic maps of the interval with a single discontinuity which involves infinite sequences of bifurcations, and is robust within the class of maps considered. This anharmonic route is one (or more accurately a countable set) out of an uncountable number of possibilities for piecewise monotonic maps \cite{G2014}, and it is natural to ask whether there are new transitions that are robust to perturbations of the system within appropriate classes of maps

In this paper we describe the boundary of the set of parameters with positive topological entropy, a set we call the boundary of chaos, for a class of one-dimensional maps with two piecewise monotonic components and a pair of plateaus. We show that 
the codimension one transitions across this boundary are either generalisations of the heteroclinic transition of circle maps \cite{MT1986} or the anharmonic transition \cite{G1993}. The construction of the boundary of chaos also shows implicitly how and where each possibility arises.

It turns out that questions about the boundary of chaos for maps with plateaus are closely connected with the structure of survivor sets of open maps. Given a map $f:\mathbb{R}\to\mathbb{R}$ and an open set $H\subset \mathbb{R}$, the open map $f_H$ is the map which coincides with $f$ on $\mathbb{R}\backslash H$ and is undefined on $H$. The set $H$ is called the hole or holes. A natural problem in the study of open maps is to characterise the dynamics of $f_H$ which avoids $H$, i.e. what can be said about the \emph{survivor set}, the points with $f_H^k(x)\notin H$, $k=0,1,2,\dots$? A related question is how this set changes as $H$ changes. In such problems, if $H$ has a simple parameterization or characterization and the dynamics of $f$ is sufficiently rich,   
then the survivor set may be large (positive Hausdorff dimension for example) or small (zero Hausdorff dimension).
What is the boundary between these regions? More detailed information about the boundary may be available, for example is the survivor set infinite or finite on the boundary, and how does this change on different parts of the boundary. This boundary is by no means the only interesting object of study. The speed at which solutions are absorbed by the hole and structure of invariant measures can have interesting properties \cite{Bahsoun2010,Bon2022,Bose2014,Bun2011,PY1979} . 

A natural example, in the sense that the dynamics is both rich and well-understood, is the doubling map, $T_2(x)=2x~(\textrm{mod}~1)$. In this case, if the hole is $(a,b)$, $a \le \frac{1}{2}\le b$, then a complete description of the boundary in the $(a,b)$ parameter space between those parameters at which the survivor set has zero Hausdorff dimension and those with zero Hausdorff dimension can be given (\cite{GS2015}, see also Corollary~\ref{cor:++} in section~\ref{sect:bdytypeA}). On the codimension one line segments of the boundary for this class of open maps, the maps have a finite number of periodic points. One of the aims of this paper is to show how this changes for the more general class of maps which cover the unit interval twice but which may have decreasing branches. We show that a new codimension one component of the boundary appears on which the survivor set has infinitely many periodic points. This is the anharmonic cascade of \cite{G1992,G1993}.     


\begin{figure}
\centering
\includegraphics[width=12cm]{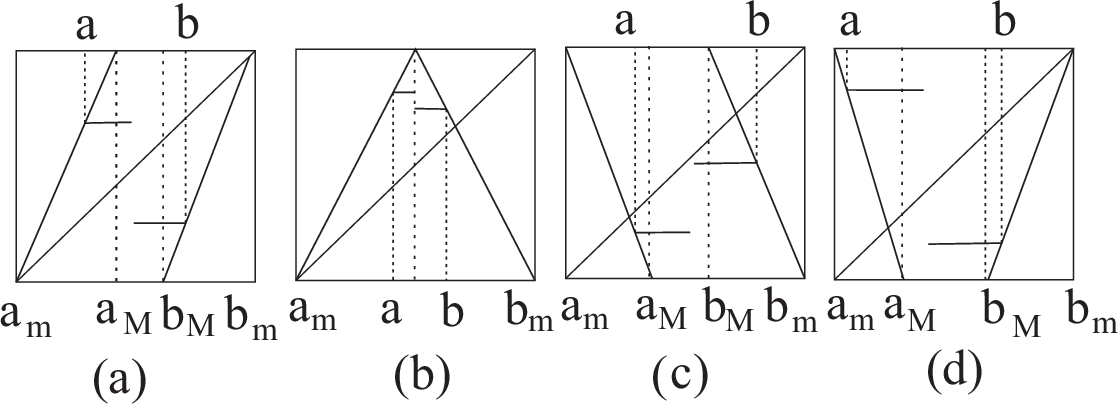}
\caption{Maps with holes and plateau maps. (a) Type A; (b) Type B; (c) Type C; (d) Type D. }
\label{fig:one}
\end{figure}


Piecewise smooth maps have a long history, partly due to their appearance as models of homoclinic phenomena \cite{GW1979,Will1979}, partly because they act as simple models of smooth dynamics with greater control over hyperbolicity \cite{Lozi1978,M1980} and partly because of their applications in engineering and biology. The transition to chaos in piecewise smooth maps with two increasing branches was described in \cite{GPTT} using ideas from the analysis of circle maps developed in \cite{MT1986}. The generalization to two monotonic branches was sketched by \cite{GLT1990} and described in more detail in \cite{G1993,G2014}. Although the motivation was from families of maps such as 
\begin{equation}
x\to \left\{ \begin{array}{ll}\mu +ax^2 & {\rm if}~x<0\\ \mu +bx^2 & {\rm if}~x>0\end{array}\right.  ,
\label{eq:smoothcases}\end{equation}
with $a,b\in \{-1,+1\}$, all of these papers are based on symbolic descriptions of orbits, with no formal proof that all symbolic possibilities are actually realized in the motivating families. Indeed, not all topologically possible dynamics do occur in these families: there are problems with homtervals (wandering intervals on which every iterate of the map is a homeomorphism) \cite{BM1991,KStP2001}. 

Open maps can be associated with maps having plateaus as shown in Figure~\ref{fig:one}. One interpretation of the results of \cite{GS2015} is that the plateau maps associated with holes for $T_2$ are full, in the sense that all possible topological dynamics compatible with a piecewise smooth map having two increasing branches, is realised. This idea is made more precise in section~\ref{sect:basic}. One of the results below (Theorem~\ref{thm:full}) is that the same is true for expanding maps with two monotonic branches whatever their orientation. 

One final historical comparison is worth making. The existence of full families, the monotonicity of kneading invariants and other related results have been proved for trapezoidal maps \cite{BMT1999} and in many ways the plateau maps defined here are the piecewise continuous analogues of these maps. From this point of view our results fit naturally into the context of simplified families of maps which make some notoriously hard problems in smooth realizations much easier to approach.         

In section~\ref{sect:basic} we describe the classes of maps studies and introduce the ideas used to describe the dynamics of these maps. This leads to the first major result (Theorem~\ref{thm:full}) that the classes of plateau map are full. In section~\ref{sect:chaos} we review the codimension one transitions to chaos in piecewise continuous maps based on the results of MacKay and Tresser \cite{MT1986} and Glendinning \cite{G1992,G1993}. In section~\ref{sect:bdycon} the boundary is shown to be connected and in section~\ref{sect:bdybdy} we describe how the boundary of chaos intersects the boundaries of the parameter space on which the maps are defined. The first steps of an inductive process to describe the boundary of chaos for plateau maps in class $C$ and class $B$ are given in sections~\ref{sect:firstC} and \ref{sect:firstB} respectively. In section~\ref{sect:bdytypeA} we describe the boundary of chaos for class $A$ maps following the results of \cite{GPTT,GS2015}, and this makes it possible to prove the main result characterising the boundary of chaos and the boundary of positive Hausdorff dimension, Theorem~\ref{thm:main}, in section~\ref{sect:mainresult}. The paper ends with the short discussion of section~\ref{sect:conclusion}. 


\section{Basic definitions}\label{sect:basic}
Throughout this section and following sections we will take limits from above and below using the notation
\[
f(x_-)=\lim_{y\uparrow x}f(y), \quad f(x_+)=\lim_{y\downarrow x}f(y),
\]
provided the limits exist.

\begin{definition}$f:\mathbb{R}\to \mathbb{R}$ is a monotonic single discontinuity map (MSDM) on $[a_m,b_m]$ if there exists $a_m< c< b_m$ such that $f|_{(a_m,c)}$ and $f|_{(c,b_m)}$ is continuous and strictly monotonic with $f(a_{m})$, $f(c_-)$, $f(c_+)$ and $f(b_{m})$ all in $[a_m,b_m]$. 
\end{definition}

Note that the value of $f$ at the point of discontinuity $c$ is left unspecified. The natural candidates are to leave it undefined, or to use continuity from either the left or the right. The definitions used below give results independent of this choice (by the strict monotonicity of the branches) and though there will be some technical complications when considering maps with plateaus. 

We will say a $MSDM$ is in class A if both branches are increasing and $f(a_m)=a_m$, $f(b_m)=b_m$; class B if the left branch is increasing and the right decreasing and $f(a_m)=f(b_m)=a_m$, class D if the left branch is decreasing and the right branch is increasing with $f(a_m)=f(b_m)=b_m$ and class C if $f$ both branches are decreasing with  $f(b_m)=a_m$ and $f(a_m)=b_m$.

The analysis of piecewise monotonic maps in one dimension uses kneading theory, see \cite{MT1988,dMvS} for details. The address, $a(x)$, of a point $x\in [a_m,b_m]$ is $0$ if $x\in [a_m,c )$ and $1$ if $x\in (c,b_m]$.  Suppose that $f^k(x)\ne c$, for all $k \ge 0$, then the kneading sequence of $x$ is
\[
k(x)=a(x)a(f(x))a(f^2(x))\dots \in \{0,1\}^{\mathbb{N}_0}.
\]
Preimages of $c$ have two associated kneading sequences, the upper and lower kneading sequences, $k (x_\pm)$, obtained by approaching $x$ from above or below through points which are not preimages of $a_M$ (this is always possible by the strict monotonicity of $f$).

We now define standard orders on sequences which respect the order of points and their iterates \cite{MT1988}. The parity operator on sequences has $p(0)=1$ if $f$ is increasing in $(a_m,c)$ and $p(0)=-1$ if $f$ is decreasing. Similarly $p(1)=1$ if $f$ is increasing on $(c,b_m)$ and $p(1)=-1$ otherwise. Define $<_s$ by 
\[
a_0\dots a_{n-1}a_n\dots <_s a_0\dots a_{n-1}b_n\dots, \quad a_n\ne b_n,
\] 
if and only if either $a_0<b_0$ or
\[\begin{array}{lll}
\prod_0^{n-1} p(a_k(x)) >0, & \textrm{and} & a_n=0, \ \ b_n=1;\ \ \textrm{or}\\
\prod_0^{n-1} p(a_k(x)) <0, & \textrm{and} & a_n=1, \ \ b_n=0.
\end{array}
\]
This order is chosen so that $x<y$ implies that $k(x_+)\le k(y_-)$. The kneading invariant of a MSDM $f$ is the pair $k_f=(k(c_-),k(c_+))$, and topological behaviour is essentially defined by this pair. The dynamics of a point $x$ under $f$ is reflected in the dynamics of the kneading sequences of $x$ under the shift map $\sigma$: $k(f(x))=\sigma k(x)$  for points that are not preimages of $c$, whilst for points that are preimages of $x$ a similar relation holds but limits from above may become
limits from below and vice versa. For consistency
\begin{equation}\label{eq:kneadconsist}
\textrm{either} \quad \sigma^rk_\pm\le_s k_-\quad \textrm{or} \quad \sigma^rk_\pm\ge_s k_+ ,
\end{equation}
for all $r\ge 0$. Any pair of sequences $(k_-,k_+)$ which satisfy (\ref{eq:kneadconsist}) and $k_-=0\dots$ and $k_+=1\dots$ is called a compatible pair.

\begin{definition}\label{def:expdc}Given $a_m<a_M \le b_M<b_m$, we will say that $f:[a_m,a_M]\cup [b_M,b_m]\to \mathbb{R}$ is an expanding double cover if
$f((a_m,a_M ))=f((b_M,b_m))=(a_m,b_m)$, $f$ is continuously differentiable on $(a_m,a_M )$ and $(b_M,b_m)$, and there exist constants $\lambda$ and $\nu$, $1<\lambda\le \nu <\infty$ such that $\lambda \le |f^\prime (x)|\le \nu$ for all $x\in (a_m,a_M)\cup (b_M,b_m)$.
\end{definition}

The slope condition together with the continuity of the derivative implies that an expanding double cover is strictly monotonic on $(a_m,a_M)$ and $(b_m,b_m)$ and we will continue to use classes A, B, C and D to describe the different types of double cover in the obvious way. We will not be interested in $f$ between $a_M$ and $b_M$ so we can consider it to be undefined in this open (and possibly empty) interval.

\begin{lemma}Suppose that $f$ is an expanding double cover. If $a_M<b_M$ then for each sequence $s\in\{0,1\}^{\mathbb{N}_0}$ there exists $x\in [a_m,b_m]$ such that $k(x)=s$. If $a_M =b_M$ then there exists $x\in [a_m,b_m]$ such that either $k(x_p)=s$ for some $p\in\{-,+\}$.
\label{lem:allseq}\end{lemma}

The proof is an exercise in the application of covering intervals (e.g. \cite{BGMY,G1994}).

\begin{definition}Suppose that $f$ is an expanding double cover. The address $\alpha$ of a symbol sequence $s\in\{0,1\}^{\mathbb{N}}$ is defined by $\alpha (s)=x$ iff $k (x_p )=s$ for some $p\in\{-,+\}$.
\end{definition}

Note that $x$ is well-defined since the cover map is expanding for every $s\in\{0,1\}^{\mathbb{N}}$. The following lemma is an immediate consequence of the definitions and construction.

\begin{lemma}Suppose that $f$ is an expanding double cover and $a\in [a_m, a_M]$ and  $b\in [b_M, b_m]$. If $a\le a_M\le b_M\le b$ with strict inequality in at least one case then the map $f_{a,b}$ defined by $f_{a,b}(x)=f(x)$ if $x\in [a_m,a)\cup (b,b_m]$, $f(a)=f(a_-)$, $f(b)=f(b_+)$ and $f_{a,b}$ is undefined on $(a,b)$ is an open map with hole $H=(a,b)$.

If $a=a_M=b_m=b$ then $f_{a,b}$ defined by $f_{a,b}(x)=f(x)$ if $x\in[a_m,b_m]\backslash\{a_M\}$ is a MSDM with discontinuity at $a_M$. In this case the topological entropy of $f_{a,b}$ is $\log~2$ and the Hausdorff dimension of the points that are not preimages of the discontinuity is one.
\label{lem:openfab}\end{lemma}

Recall that the survivor set $S$ of a map with holes, $f_{a,b}$ is   
\begin{equation}\label{eq:S}
S=\left\{x\in X~|~f_{a,b}^n(x)\notin (a,b)~\textrm{ for~all }~n\geq0\right\}.
\end{equation}

\begin{definition}A family of maps $f_{a,b}$ defined via the construction of Lemma~\ref{lem:openfab} is a complete family if all possible $(a,b)$ occur, i.e. it is the family
\begin{equation}\label{eq:complete}
 \{ f_{a,b}~|~(a,b)\in [a_m,a_M]\times [b_M,b_m]\}.
\end{equation}\end{definition}

\begin{definition}\label{def:plateau}
Fix $c\in [a_M ,b_M]$. Then a family of plateau maps is $F_{a,b}$ where $F_{a,b}(x)=f_{a,b}(x)$ except in the hole where $F_{a,b}(x)=f_{a,b}(a)$ if $x\in (a,c)$ and   $F_{a,b}(x)=f_{a,b}(b)$  if $x\in (c,b)$ provided these are defined. $\{F_{a,b}~|~(a,b)\in M\}$ is a complete family if $M= [a_m,a_M ]\times [b_M, b_m]$.
\end{definition}

Note that this leaves $F_{a,b}$ undefined at $x=c$. If the upper and lower kneading sequences are defined then this pair is the kneading invariant of $F_{a,b}$ as with cases with strictly monotonic branches. If this is not possible, i.e. if the limits as $x$ tends to $a$ from below or $b$ from above do not exist, then we make the arbitrary choice $F_{a,b}(c)=F_{a,b}(c_-)$ so that every plateau map has a kneading invariant. This does not change the space of compatible kneading invariants, but it does alter the details of the relationship between the possible dynamics and the corresponding kneading invariants. Since this level 
of detail is not used in this paper this simple solution for the definition of kneading invariants is sufficient for our classification.   

The connection between the study of open maps and their corresponding plateau map is described by a classic relation between the Hausdorff dimension of the survivor set of the open map, $\mathcal{H}_D(S)$, and the topological entropy of the plateau map, $h_{top}(F_{a,b})$. 

\begin{proposition}
Let $f$ be an expanding double cover and for $(a,b)\in [a_m,a_M]\times [b_M,b_m]$ let $f_{a,b}$ and $F_{a,b}$ be corresponding open maps (Lemma~\ref{lem:openfab}) and plateau maps (definition~\ref{def:plateau}). If the survivor set of $f_{a,b}$ is $S$ then
\[
\mathcal{H}_D(S)>0\quad {\rm if\ and \ only\ if}\quad h_{top}(F_{a,b})>0.
\]
\end{proposition}
   
The proof of this proposition for invariant sets of the doubling map, $T_2$, can be found in Furstenberg \cite[Proposition III.1, p.~36]{furstenberg1967disjointness}. Furstenberg's proof extends to the other expanding double covers with little difficulty, see \cite[Section 1.2.2, p.~39]{Hege2022} for details of the case $\lambda=\nu =2$ in Definition~\ref{def:expdc}. The more
general case follows the same argument, using the fact that the length $|J|$ of an interval $J$ on which $f^n$ is a homeomorphism and $f^n(J)=(a_m,b_m)$ is bounded above and below by $(b_m-a_m)\lambda^{-n}$ and $(b_m-a_m)\nu^{-n}$ respectively.

\begin{definition}Given a parameter space $M$, a family of MSDMs $\{f_\mu~|~\mu\in M\}$ or plateau maps $\{F_{a,b}~|~(a,b)\in M\}$ is called a full family in class T (T=A,B,C,D) if for every compatible pair $(k_-,k_+)$ there exists $m\in M$ such that the kneading invariant of $f_m$ (for MSDMs) or $F_m$ (for plateau maps) is $(k_-,k_+)$.\end{definition}

Thus full families can exhibit `every' type of topological dynamics consistent with the class of the map. Nice classes of maps, e.g. those used in \cite{GPTT,G1993, G2014, GLT1990} do not have this property \cite{KStP2001}. Theorem~\ref{thm:full}, at the end of this section, shows that it is possible to create families of plateau maps which are full. It also shows that the choice of definition of a plateau map at $c$ is not important for the main results here.

The next result is the analogue of the equivalent result for trapezium maps in \cite{BMT1999}.

\begin{theorem}Complete families of plateau maps are full families.\label{thm:full}\end{theorem}

\emph{Proof:} Choose any compatible pair $(k_-,k_+)$ in the appropriate class and let $f$ be a double cover of the same class. The points $a_m$ with $k(a_{m-})=k_-$ and $b_m$ with $k(b_{m+})=k_+$ exist by Lemma~\ref{lem:allseq} and each avoids the plateau by the definition of compatible sequences, (\ref{eq:kneadconsist}). So the kneading invariant of $F_{a,b}$ is $(k_-,k_+)$.
\newline\rightline{$\square$}


\section{Codimension one transitions to chaos}\label{sect:chaos}
There are three known codimension one transitions to chaos (in the sense of positive topological entropy) for piecewise monotonic maps of the interval. The heteroclinic transition and the circle intermittency transition \cite{GPTT,MT1986} both derive drom the study of non-invertible maps of the circle, and the anharmonic route \cite{G1992,G1993} requires at least one of the branches to be decreasing. One of the main results of this paper is that the only codimension one transitions in full families of plateau maps are the heteroclinic transition and the anharmonic transition. The heteroclinic transition is an abrupt transition in that there is a finite number of periodic points for maps on the boundary of chaos, and a small perturbation leads to maps with infinitely many periodic orbits. The anharmonic route \cite{G1992,G1993} involves an infinite cascade of bifurcations so that on the boundary there is a countably infinite number of periodic orbits. In some ways this is the piecewise smooth equivalent of period-doubling for continuous maps (which is codimension two in the space of piecewise smooth maps we are considering here), and in the simplest case the cascade is a sequence of doubling the period and subtracting one followed by doubling the period and adding one. A heuristic explanation for the absence of the circle intermittency route as a codimension one phenomenon for the plateau maps considered here is given at the end of section~\ref{sect:bdytypeA}. 

Although these are both fundamental codimension one transitions to chaos they are, perhaps, not as broadly known as they deserve to be, so we provide the statements here along with a sketch proof of the transitions. For the anharmonic route we give a new proof independent of the kneading theory developed in \cite{G1993}.

The heteroclinic transition is based on the idea of well-ordered periodic orbits. These are periodic orbits that have the same symbolic description as those arising in circle maps, and they do not imply the existence of other periodic orbits. 

\begin{definition}\label{def:balanced}A word $\omega\in \{0,1\}^\infty$ is balanced if given any two segments of the same length the number of 1s in each segment differs by at most one. A periodic orbit is well-ordered if its kneading sequence is balanced.
\end{definition}

Non-periodic balanced words are also called Sturmian \cite{V2003}. For $0<p<q$ with $p$ and $q$ coprime let $\omega_{p/q}$ be a periodic balanced word of period $q$ with $p$ 1s in each period. Let $[\omega ]_n$ denote the first $n$ symbols of $\omega$ and using the lexicographical order define  
\[\begin{array}{lll}
r^-_{p/q}=[0M_{p/q}]_q, &{\rm where} & M_{p/q}=\max_k\sigma^k\omega_{p/q},\\
r^+_{p/q}=[1m_{p/q}]_q, & {\rm where} & m_{p/q}=\min_k\sigma^k\omega_{p/q}.
\end{array}\] 

\begin{theorem}(MacKay-Tresser) Suppose $f_\mu$ is a class A continuous family of maps and that if $\mu =0$ then $p_1<  \dots <p_r<c<p_{r+1}<\dots <p_q$ are points on a well-ordered unstable periodic orbit with
\begin{equation}
f_\mu^q(c_-)=p_{r+1}, \quad  f_\mu^q(c_+)\in (0,p_{r+1}), \quad \frac{\, d}{d\mu} [f_\mu (c_-)-p_{r+1}]{\big|}_{\mu =0}=b\ne 0.
\end{equation}
If $b>0$ then there is a neighbourhood of $\mu =0$ such that $f_\mu$ has positive topological entropy if $\mu>0$ and zero topological entropy if $\mu <0$ in this neighbourhood. If $b<0$ then the roles of $\mu<0$ and $\mu>0$ are interchanged.
\label{thm:hettrans}
\end{theorem}

Of course, there is an equivalent statement if
\begin{equation}
f_\mu^q(c_+)=p_{r}, \quad  f_\mu^q(c_-)\in (p_r,0), \quad \textrm{and}\quad \frac{\, d}{d\mu} [f_\mu (c_+)-p_r]{\big|}_{\mu =0}\ne 0.
\end{equation}
using the correspondence obtained by reversing the orientation of the interval. Note that the second and third conditions are open and so if they are satisfied then they are satisfied for nearby families, and so there is only one condition to be satisfied and solutions persist (by the implicit function theorem for smooth enough families) and we see that the transition is indeed codimension one. Note also that the second condition cannot be satisfied by expanding maps, so the class of maps in which the transition occurs is not 'all' such maps (the same is true of period-doubling in the continuous case).

\emph{Sketch proof of Theorem~\ref{thm:hettrans}}

\begin{figure}
\centering
\includegraphics[width=10cm]{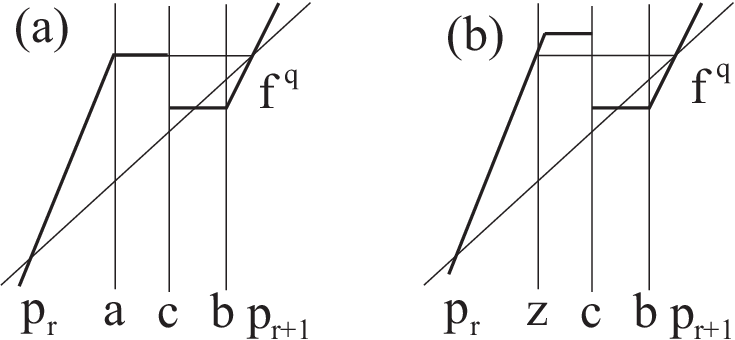}
\caption{The heteroclinic transition illustrated for a family of plateau maps. (a) On the boundary with $f^q(c_-)=p_{r+1}$: all orbits map either to the stable fixed point in the plateau or to the unstable fixed point $p_{r+1}$ under $f^q$; (b) immediately after the transition, with $z<a<c$ such that $f^q(z)=p_{r+1}$. }
\label{fig:het}
\end{figure}

The induced map $f_0^q$ restricted to $[p_r,p_{r+1}]$ or indeed any of the intervals defined by successive points of the $(q-r,q)$-orbit which are topologically conjugate to this central region, takes the form shown in Figure~\ref{fig:het} and there is a stable periodic orbit in $(c,p_{r+1})$ and there may be other such stable periodic orbits with the same symbolic description. Assume that the well-ordered orbit has been chosen so that there are no stable periodic orbits in $(p_r,c)$, and note that hyperbolicity implies that the orbit persists for all sufficiently small perturbations from $\mu =0$. 


If $\mu<0$ with $|\mu|$ small then $f_\mu^q(c_-)\in (c,p_{r+1})$ and the dynamics is essentially unchanged from the boundary case $\mu =0$.

For any $\mu >0$ and $|\mu |$ small, there exists $z\in (p_r,c)$ such that $f_\mu^q(z)=p_{r+1}$ and so
\begin{equation}\label{eq:coverleft}
f_\mu^q([p_r,z])=[p_r, p_{r+1}]\supset [p_r,c].
\end{equation}
Note that preimages of $c$ under $f_\mu^q$ accumulate on the right hand side of $p_r$ and hence on the right hand side of $p_s$ for any $s$. Now, $f_\mu^q([z,c))=[p_{r+1}, u)$ for some $u>p_{r+1}$, and since, as just shown, preimages of $c$ accumulate on $p_{r+1}$ from above, we can choose $v$ such that $p_{r+1}<v<u$ such that there exists $m>0$ with $f_\mu^m$ monotonic on $(p_{r+1},v)$, $f_\mu^m(p_{r+1})=p_r$ and $f_\mu^m(v)=c$. Thus $f^{q+m}_\mu$ is continuous on $(z,c)$ and 
\begin{equation}\label{eq:coverright}
f_\mu^{q+m}([z,c)])\supset [p_r,c].
\end{equation}
Equations (\ref{eq:coverleft}) and (\ref{eq:coverright}) imply that there is a topological horseshoe in the map and hence that the entropy is positive.
\newline\rightline{$\square$}

The anharmonic route was described in \cite{G1992} and precise statements together with the detailed proof can be found in \cite{G1993}. It applies to maps in class $B$, $C$ and $D$, and uses induced maps to show that there can be robust, i.e. codimension one, transitions to chaos involving infinitely many bifurcations. An induced map argument is used and since this echoes the construction used in the two parameter problems here it is worth restating the basic argument in terms of plateau maps. This has the advantage that it is unnecessary to go through the detailed kneading theory used in \cite{G1993}.

\begin{definition}Suppose that $f:[a_m,a_M]\cup [b_M,b_m]\to \mathbb{R}$ is an expanding double cover. For fixed $b\in [b_M,b_m]$ a family of plateau maps $\{F_{a,b}~|~a\in M\}$ as defined in Definition~\ref{def:plateau} is $a$-full if $M=[a_m,a_M]$. \label{def:afull}\end{definition}   

\begin{lemma}Suppose that $\{F_{a,b}\}$ is an $a$-full family of plateau maps of class $C$ and $b\in [b_M,x^*]$, where $x^*$ is the fixed point of the map $(b_M,b_m)$. Then there exist $a_1$, $a_2$, with $a_m<a_1<a_2<a_M$ such that if $a\in [a_1,a_2]$ the induced map
\begin{equation}
R_CF_{a,b}=\left\{\begin{array}{ll}F_{a,b}^2(x) &{\rm if}\ x\in [a_1,c)\\ F_{a,b}(x) &{\rm if}\ x\in (c,F(a_2)] \end{array}\right. ,
\label{eq:inducedCB}\end{equation}  
is an $a$-full family of plateau maps in class B with plateau $(a,b)$. Moreover, if $a\in[a_m,a_1)$ then $h_{top}(F_{a,b}=0$ and if $a\in (a_2, a_M]$ then $h_{top}(F_{a,b}>0$.

Similarly, suppose that $\{F_{a,b}\}$ is an $a$-full family of plateau maps of class $B$ and $b\in [b_M,x^*]$. Then there exist $a_1$, $a_2$, with $a_m<a_1<a_2<a_M$ such that if $a\in [a_1,a_2]$ the induced map
\begin{equation}
R_BF_{a,b}=\left\{\begin{array}{ll}F^3_{a,b}(x) &{\rm if}\ x\in [a_1,c)\\ F_{a,b}(x) &{\rm if}\ x\in (c,F(a_1)] \end{array}\right. ,
\label{eq:inducedBC}\end{equation}  
is an $a$-full family of plateau maps in class C with plateau $(a,b)$.
\label{lem:anharminduced}\end{lemma}

\begin{figure}
\centering
\includegraphics[width=10cm]{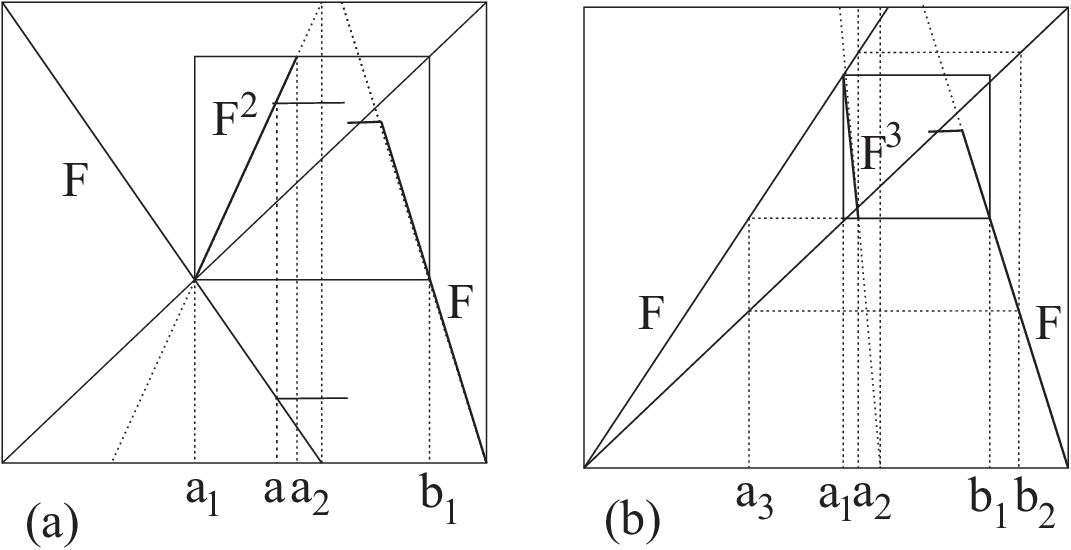}
\caption{Construction of induced maps. (a) (\ref{eq:inducedCB}); and (b) (\ref{eq:inducedBC}). }
\label{fig:anharm}
\end{figure}

\emph{Proof:} Figure~\ref{fig:anharm} shows the construction of the two families, and to some extent no more explanation is required. However, let us make this more explicit. We will drop the subscripts $a,b$ when the choice of map is clear. 

Suppose that $\{F_{a,b}\}$ is an $a$-full family of plateau maps of class $C$ and $b\in [b_M,x^*]$, where $F(x^*)=x^*$, $x*>b_M$ ($x^*$ exists by the intermediate value theorem and the relationships $F(b_M)=b_m$ and $F(b_m)=a_m$ in the definition of class $C$ maps). Since $F(a_m)=b_m$ and $F(a_M)=a_m$ there is a fixed point $a_1\in (a_m,a_M)$. Let $b_1>b_M$ be the point such that $F(b_1)=a_M$ and $a_2<a_M$ be the point such that $F(a_2)=b_1$. Then $F^2|\, [a_1,a_2]$ and $F|\, [b_M,b_1]$ is an expanding double cover and the family of plateau maps (\ref{eq:inducedCB}) is an $a$-full family of maps in class $B$. The plateau regions of the induced map and $F_{a,b}$ coincide since $F([a_1,a_M])\subseteq [a_m,a_1]$ and $F|\, [a_m,a_1]$ is a homeomorphism if $a\ge a_1$. The argument for the claims about the topological entropy follow standard arguments as in the sketch proof of Theorem~\ref{thm:hettrans} and is omitted, see \cite{G1993} for details.


In the case that $\{F_{a,b}\}$ is an $a$-full family of plateau maps of class $B$ and $b\in [b_M,x^*]$, the argument is similar in spirit, but the construction involves more points.
As in the previous case define $x_1>b_M$ such that $F(x_1)=a_M$ and $x_2<a_M$ such that $F(x_2)=a_M$. Note that there exists $a_1\in (a_m,x_2)$ such that $b_1=F(a_1)>x_1$ and $F(b_1)=a_1$ and $x_2<a_1$ since $F$ is increasing on $[a_m,a_M]$. 

Now define $x_3>x_1>b_M$ such that $F(x_3)=x_2$ and $x_4<a_M$ such that $F(x_4)=x_3$. Thus $F^3(x_4)=a_M$ with $F(x_4)>x_1>b_M$ and  $F^2(x_4)=x_2<a_1$. By the definition of class $B$ maps, $F(a_M)=b_m$, $F^2(a_M)=F^3(a_M)=a_m$. Thus $F^3|\, [x_4,a_m]$ is a decreasing homeomorphism with $F^3(x_4)=a_M$ and $F^3(a_M)=a_m$. Hence there exists $a_2\in (x_4,a_M)$ such that $F^3(a_2)=a_1$ and so $F^3([a_1,a_2])=[a_1,b_1]$.

The argument is now as in the previous case: the induced map (\ref{eq:inducedBC}) is an $a$-full map in class $C$ with $a\in [a_1,a_2]$ and the entropy statement again follows the same style of argument as in the proof of Theorem~\ref{thm:hettrans}, see \cite{G1993} for details.
\newline\rightline{$\square$}

Just as in the case of period-doubling cascades, the operators $R_B$ and $R_C$ can be associated with replacement rules on symbol sequences. Suppose that $F_{a,b}$ is an $a$-full family of plateau maps of class $C$ with $b\in (b_M,x^*)$ as in Lemma~\ref{lem:anharminduced}. Then $R_BR_CF_{a,b}$ is also an $a$-full family of plateau maps of class $C$ with $b\in (b_M,x^*)$, and if the kneading sequences in $\{0,1\}^{\mathbb{N}_0}$ of a point $x$ under $R_BR_CF_{a,b}$ is known, then the kneading sequence of $x$ under $F_{a,b}$ is obtained using the replacement rule
\begin{equation}
(0,1)\to (00100,1) ,
\label{eq:replaceC}\end{equation}
which is $(00,1)\to (00,1)$ for $R_C$ followed by $(0,1)\to (010,1)$ for $R_B$.
Similarly for case $B$ in Lemma~\ref{lem:anharminduced} the replacement rule is obtained by reversing the order of these two replacement operations:
\begin{equation}
(0,1)\to (010010,1) .
\label{eq:replaceB}\end{equation}

\begin{theorem}Suppose that $\{F_{a,b}\}$ is an $a$-full family of plateau maps of class $C$ and $b\in [b_M,x^*]$, where $F(x^*)=x^*$. Then there exists a unique $a_C\in [a_m,a_M]$ such that $h_{top}(F_{a,b})=0$ if $a\le a_C$ and $h_{top}(F_{a,b})>0$ if $a> a_C$. If $a=a_C$ then the kneading sequence of $F_{a,b}$ is $(k_+,k_-)$ with $k_+=1^\infty$ and $k_-$ is the sequence obtained by repeated application of the replacement rule (\ref{eq:replaceC}). Moreover, if $a=a_C$ then $F_{a,b}$ has a countable number of periodic orbits with periods $p_n$, $n=0,1,2,\dots $ , where  
\begin{equation}
p_n=\frac{1}{3}(4.2^n-(-1)^n).
\label{eq:pnC}\end{equation}
Similarly, if $\{F_{a,b}\}$ is an $a$-full family of plateau maps of class $B$ and $b\in [b_M,x^*]$, where $F(x^*)=x^*$ then there exists a unique $a_B\in [a_m,a_M]$ such that $h_{top}(F_{a,b})=0$ if $a\le a_B$ and $h_{top}(F_{a,b})>0$ if $a> a_B$. If $a=a_B$ then the kneading sequence of $F_{a,b}$ is $(k_+,k_-)$ with $k_+=1^\infty$ and $k_-$ is the sequence obtained by repeated application of the replacement rule (\ref{eq:replaceB}). Moreover, if $a=a_B$ then $F_{a,b}$ has a countable number of periodic orbits with periods $p_n$, $n=0,1,2,\dots $ , where  
\begin{equation}
p_n=\frac{1}{3}(5.2^n+(-1)^n).
\label{eq:pnB}\end{equation}
\label{thm:anharm}\end{theorem}
  
\emph{Proof:} In the case of class $C$, induction on Lemma~\ref{lem:anharminduced} proves that the induced maps $(R_BR_C)^nF_{a,b}$ are $a$-full families of class $C$ maps with $\in (b_M,x^*)$ on nested intervals $[a_1^{(n)},a_2^{(n)}]$ and that the map has zero entropy below $a_1^{(n)}$ and positive entropy above $a_2^{(n)}$. The set
\[ a_C=\bigcap_{n\ge 0}[a_1^{(n)},a_2^{(n)}]
\]
is clearly non-empty and by definition the double covers are expanding so by considering the slopes of the induced maps $a_2^{(n)}-a_1^{(n)}$ decreases faster than exponentially, hence $a_C$ is a single point.  

At $a=a_C$ the induced maps $(R_BR_C)^nF_{a,b}$ can be defined for all $n>0$ and so the kneading invariant $k_-$ is obtained by applying (\ref{eq:replaceC}) infinitely often. 

By definition, $F_{a,b}$ has a boundary orbit of period two, which we take as given. To apply $R_C$, the fixed point in $x<a_M$ must exist which has symbolic description $(0)^\infty$ and period $p_0=1$. To apply $R_B$ to $R_CF_{a,b}$ the period two orbit of $R_CF_{a,b}$ $(01)^\infty$ must exist, which using the replacement $(0,1)\to (00,1)$ is period 3 ($p_1=3$) with kneading sequence $(001)^\infty$. The left hand branch of $(R_BR_C)F_{a,b}$ is obtained using (\ref{eq:replaceC}) and so the fixed point on this branch, which allows the induced map $R_C(R_BR_C)F_{a,b}$ to be defined has period 5 ($p_2=5$) with kneading sequence $(00100)^\infty$. 

Mote generally suppose $p_n$ is the period under $F_{a,b}$ of the fixed point for $(R_BR_C)^nF_{a,b}$ (so $n$ is even). Then the period of the period  two orbit of $R_C(R_BR_C)^nF_{a,b}$ which makes it possible to define $R_BR_C(R_BR_C)^nF_{a,b}$ has kneading sequence $(001)^\infty$ for $(R_BR_C)^nF_{a,b}$, i.e period $p_{n+1}=2p_n+1$ for $F_{a,b}$. Finally, the fixed point of 
$(R_BR_C)^{n+1}F_{a,b}$ has kneading sequence $(00100)^\infty$ for $(R_BR_C)^nF_{a,b}$, i.e. 
\[
p_{n+2}=4p_n+1=2p_{n+1}-1.
\]
Since $n$ is even, the expressions for $p_{n+1}$ and $p_{n+2}$ can be written as
\begin{equation}
p_{n+1}=2p_n-(-1)^n, 
\label{eq:2npm}\end{equation} 
which is valid for $n$ even and $n$ odd. With the initial condition $p_0=1$ a straightforward calculation gives (\ref{eq:pnC}).

The entropy statements follow from Lemma~\ref{lem:anharminduced} and the case class $B$ maps uses precisely the same style of argument.
\newline\rightline{$\square$}

The sequences (\ref{eq:pnC}) and (\ref{eq:pnB}) are not the only sequences of periods associated with this transition since
the robust renormalization can associated with an anharmonic cascade can arise after a different set of renormalizations.

         
\section{The boundary is connected}\label{sect:bdycon}     
Given a family of maps $\{f_\mu ~|~ \mu\in M\}$ where $M$ is a connected open subset of $\mathbb{R}^n$, we are interested in the way $M$ divides into regions with positive topological entropy and zero topological entropy. Equivalently, given a family of open maps the issue is to characterize the boundary between regions in which the survivor set has positive Hausdorff dimension and those in which the survivor set has Hausdorff dimension equal to zero.

\begin{definition}Given a family of maps $\{f_\mu ~|~ \mu\in M\}$ where $M$ is a connected open subset of $\mathbb{R}^n$, we say $\mu_0\in M$ is on the boundary of chaos if in every open neighbourhood of $\mu_0$ there exist $\mu_1$ and $\mu_2$ such that $h_{top}(f_{\mu_1})>0$ and $h_{top}(f_{\mu_2})=0$.
\end{definition}

Clearly if this is to be a useful definition then the maps need to have some regularity with respect to the parameters, but since this is automatically the case in the examples we consider, we will not over-generalise the definition here.
     
Suppose that $\{F_{a,b}\}$ is a complete family of plateau maps. To simplify notation write $h_{top}(F_{a,b})=h_{top}(a,b)$. Clearly if $(a_1,b_1)\subseteq (a_2,b_2)$ then $h_{top}(a_1,b_1)\ge h_{top}(a_2,b_2)$.

\begin{lemma}\label{lemmaAlpha}There exists $c_1,c_2$ such that if $c_1<a+b<c_2$ then for each $c\in (c_1,c_2)$, 
\[
\alpha (c)=\max_{a+b=c}\{a~|~h_{top}(a,b)=0\}
\]
exists and $h_{top}(a,c-a))>0$ if $a>\alpha$ and $h_{top}(a,b)=0$ if $a\le \alpha$.
\end{lemma}

This follows immediately from the fact that $a_1<a_2$ implies that $c-a_1>c-a_2$ so $(a_2,c-a_2)\subset (a_1,c-a_1)$.

Note in the standard case A we may take $c_1=0.75$ so that if $b=0.5$ then $a=0.25$ and $c_2=1.25$ so that if $a=0.5$ then $b=0.75$.

\begin{lemma}\label{lem:cc}
Let $\alpha (c)$ be as above. The set
\[
G=\{a+b,\alpha(a+b)~|~a+b\in[a_m+b_M,a_M+b_m]\}. \]
is a continuous graph in the $(a+b,a)$-plane.
\end{lemma}

\emph{Proof:} Think of $G$ as a curve parametrised by $c=a+b$ so by Lemma~\ref{lemmaAlpha} each $c$ gives a unique point $(\alpha (c),\beta (c))$. Pick $c_0\in (c_1,c_2)$ and $\epsilon >0$. Then $|c-c_0|<\epsilon$ implies $|\alpha (c)-\alpha (c_0)|<\epsilon$ and $|\beta (c)-\beta (c_0)|<\epsilon$. Consider the lines $a+b=c_0$ and $a+b=c$ with $0<c-c_0<\epsilon$. If $a+b=c$ and $a>\alpha (c_0)+\epsilon$ then (as $c-a=b<\beta_0-\epsilon$), $h_{top}(a,b)>0$ so $\alpha(c)\leq\alpha(c_0)+\epsilon$. 

Similarly, if $a+b=c$ and $a<\alpha (c_0)$ then (as $c-a=b>\beta_0$), $h_{top}(\alpha,\beta )=0$ therefore $\alpha(c)\geq\alpha(c_0)$. So $\alpha$ and $\beta$ are continuous and from Lemma \ref{lemmaAlpha}, $G$ is a continuous graph on the $(a+b,a)$-plane
\newline\rightline{$\square$}

The closure of the part of $G$ which lies in the interior of the parameter box is the boundary of chaos.
 
\section{The boundary of chaos on the boundary of parameter space}\label{sect:bdybdy}
In this section we will determine the intersection of the boundary of chaos with the boundary of the parameter boxes. These results will be important when showing that the boundary links nicely through different renormalization regions in later sections.

For each class of map we define the parameter box to be the rectangle
\[
R=\{(a,b)~|~a_m\le a\le a_M, \ \ b_M\le b\le b_m\}.
\]
So $R$ includes all parameters on which the plateau map is defined and $\{F_{a,b}~|~(a,b)\in R\}$ is a complete family of plateau maps (Definition~\ref{def:plateau}) and is therefore full (Theorem~\ref{thm:full}). The boundary of each parameter box has four edges, the four sides to the rectangle. The next two results describe the topological entropy of the maps on these edges.

\begin{lemma}On the boundaries of $R$ with $a=a_m$ or $b=b_m$, $h_{top}(a,b)=0$.\label{lem:zerobdy}\end{lemma}

\emph{Proof:} If $a=a_m$ then all points in $x<a_M$ map to the bounding fixed point or point of period two in either one or two iterations. Thus any complex behaviour must lie in $x>b_M$. But $F_{a,b}$ is monotonic in $x>b_M$ so either tends to a fixed point or a point of period two. Hence the map has zero topological entropy. The argument if $b=b_m$ is similar. 
\newline\rightline{$\square$}

\begin{lemma}On the boundaries of $R$ with $a=a_M$ or $b=b_M$, 
\begin{itemize}
\item[(i)] In the case A, $h_{top}(a,b)>0$ if $b=b_M$ and $F_{a,b}(a)>b_M$ or if $a=a_M$ and $F_{a,b}(b)<a_M$. Otherwise $h_{top}(a,b)=0$.
\item[(ii)] In the case B, $h_{top}(a,b)>0$ if $b=b_M$ and $F_{a,b}(a)>b_M$ or if $a=a_M$ and $F_{a,b}(b)>a_m$. Otherwise $h_{top}(a,b)=0$.
\item[(iii)] In the case D, $h_{top}(a,b)>0$ if $b=b_M$ and $F_{a,b}(a)<b_m$ or if $a=a_M$ and $F_{a,b}(b)<a_M$. Otherwise $h_{top}(a,b)=0$.
\item[(ii)] In the case C, $h_{top}(a,b)>0$ if $b=b_M$ and either $F_{a,b}(a)<c$ or both $F_{a,b}(a)>c$ and $F_{a,b}^2(a)>b_M$; or if $a=a_M$ and $F_{a,b}(b)>c$ or both $F_{a,b}(b)<c$ and $F_{a,b}^2(b)<a_M$. Otherwise $h_{top}(a,b)=0$.
\end{itemize}
\label{lem:posbdy}\end{lemma}

\emph{Proof:} For convenience we will again drop the subscripts $(a,b)$ unless there is an ambiguity.

(i) If $b=b_M$ and $F(a)\le b_M$ then all orbits starting in $x\le a_M$ either tend to a fixed point (in the plateau or the image of the plateau). Hence as in the proof of Lemma~\ref{lem:zerobdy} the topological entropy is zero. If $F(a)>b_M$ then there exists $x<a$ such that $F(x)=b_M$ and so $F^2((x,a))=(a_m,F^2(b_M))$. Since $F$ is increasing with slope greater than one and $F(a_m)=a_m$, there exist $n>1$ and $a_1\in (a_m,F^2(a))$ such that $F(a_1) \in (a_1,F^2(a_1))$ and $F^n(a_1)=b_M$, with $F^k(a_1)<a_M$, $k=1, \dots ,n-1$. Let $I_0=(a_m,a_1)$ and $I_1=(a_1,F(a_1))$. Then $F^n(I_1)= (r,F(b_M))$ and so $I_0\cup I_1\subseteq F^{n+1}(I_1)$. By definition $I_0\cup I_1\subseteq F(I_0)$. Hence by standard techniques, $h_{top}(F)>0$. 

The case $a=a_M$ is equivalent after reversing the orientation of the interval.

\medskip
(ii) If $b=b_M$ and $F(a)\le b_M$ then all orbits starting in $x\le a_M$ either tend to a fixed point and $h_{top}(F)=0$. If $F(a)>b_M$ there exists $x<a$ such that $F(x)=b_M$ and $F^3(x,a)=(a_m,F^3(a))$. We can now argue that $h_{top}(F)>0$ using the same technique as in the proof of (i). 

If $a=a_M$ and $F(b)>a_m$ then $F((b,1))=(a_m,F(b))$ which contains preimages of $a_M$ so the same argument can be used to construct two intervals with positive entropy.

\medskip
(iii) This case is equivalent to (ii) after reversing the orientation of the interval.

\medskip
(iv) If $b=r$ then the second iterate of the double cover in $x<a_M$ has two increasing branches that each cover $[a_m,1]$. The left hand branch involves points that map to $x\ge b_M$ in one iteration, and so they are not affected by the plateau until $a$ is sufficiently small. Hence this branch together with the branch in $x>b_M$ form a type B double cover for the induced map $G(x) =F^2(x)$ in $x<c$ and $G(x)=F(x)$ in $x>c$. The condition for zero entropy of $g$ is given by part (ii) of the lemma, i.e. $G(a)\le a_m$. This is equivalent to $F(a)>c$, $F^2(a)\le b_M$. Otherwise $G$ is well defined and has positive entropy. 

If $a=a_M$ the proof is the same after a reversal of the orientation of the interval.
\newline\rightline{$\square$}

The results of Lemma~\ref{lem:posbdy} can be re-interpreted on a symbolic level.
\begin{itemize}
\item In case A the boundary of chaos intersects $a=a_M$ at $(a_M ,b_d)$ where $F(b_d)=a_M$ and the corresponding kneading invariant is $(k_-,k_+)=(01^\infty, 101^\infty )$.  Similarly, the boundary of chaos intersects $b=b_M$ at $(a_d,b_M)$ where $F(a_d)=b_M$ and the corresponding kneading invariant is $(k_-,k_+)=(010^\infty, 10^\infty )$.
\item In case B the boundary of chaos intersects the box boundary at the `corner' $(a_M,b_m)$ and $(a_d,b_M)$ where $a_d$ is defined by $F(a_d)=b_M$.  At $(a_M ,b_m)$, $(k_-,k_+)=(010^\infty,10^\infty)$ and at $(a_d,b_M)$ the kneading invariant is $(k_-,k_+)=(0110^\infty ,110^\infty )$ if $b_M>c$ or if $b_M=c$ and $a\downarrow a_d$. In the latter case $(k_-,k_+)=(0^\infty ,110^\infty )$ as $a\uparrow a_d$.
\item In case C the boundary of chaos intersects the box boundary at $(a_M,b_d)$, where $F^2(b_d)=a_M$ and $F(b_d)<c$ with $(k_-,k_+)=(00(10)^\infty ,100(01)^\infty )$  and $(a_d,b_M)$ where $F^2(a_d)=b_M$ and $F(a_d)>c$ with $(k_-,k_+)=(011(10)^\infty ,11(01)^\infty )$. 
\item Case D is equivalent to case B by symmetry.
\end{itemize}

   
\section{First renormalisation for class $C$}\label{sect:firstC}
The main tool of the next sections is the idea of induced maps. These maps will be chosen so that the entropy of the original map is positive if and only if the entropy of the induced map is positive. Hence questions about the boundary of entropy for the original map become questions about the boundary of entropy for the induced map in the parameter regions on which these induced maps are defined.

\begin{definition}\label{def:ind}
Suppose $F:[a_m,b_m]\to [a_m,b_m]$ is a complete family of plateau maps in class $T$, $T\in\{A,B,C,D\}$. If there exist $n,m>0$ with $n+m\ge 2$, and $a_g$, $b_g$ with $a_g<c<b_g$ such that $F^n(a_g)$ and $F_m(b_g)$ are in $\{a_g,b_g\}$.  Then $G:[a_g,b_g]\to [a_g,b_g]$ is an induced map if 
 \begin{equation}\label{eq:indmap}
G(x)=\left\{\begin{array}{ll}F^n(x) & \textrm{if}~x\in [a_g,c)\\F^m(x) & \textrm{if}~x\in (c,b_g]\end{array},\right.
\end{equation}
maps $[a_g,b_g]$ into itself and $F^k(x)\notin [a_g,b_g]$ if $x\in (a_g,c)$, $1\le k<n$ and if $x\in (c,b_g)$, $1\le k<m$. 
\end{definition}
The process of defining induced maps in regions of parameter space is called renormalization.

If $G$ is an induced map then
\begin{equation}\label{eq:leftover}
K=[a_m,b_m]\backslash \left(\bigcup_0^{n-1}\overline{F^k((a_g,c))}\cup \bigcup_0^{m-1}\overline{F^k((c,b_g))}\right)
\end{equation}
is a (possibly empty) finite union of intervals which are mapped over themselves by $F$ defining a finite Markov partition.

For ease of notation we will refer to an induced map $G$ defined by (\ref{eq:indmap}) as $G=(F^n,F^m)$. The regions of parameter space on which given induced maps can be defined will be denoted by $X(\omega_-,\omega_+)$, where $X$ denotes the type of the induced map (i.e. A, B, C or D), and $\omega_\pm$ are the symbol sequences of the addresses of the iterates of $(a_g,c)$ and $(c,b_g)$, starting in $[a_g,b_g]$, so their lengths $|\omega_\pm |$ denote the number of iterates the intervals take to return ($n$ and $m$ of (\ref{eq:indmap})).  

\begin{figure}
\centering
\includegraphics[width=6cm]{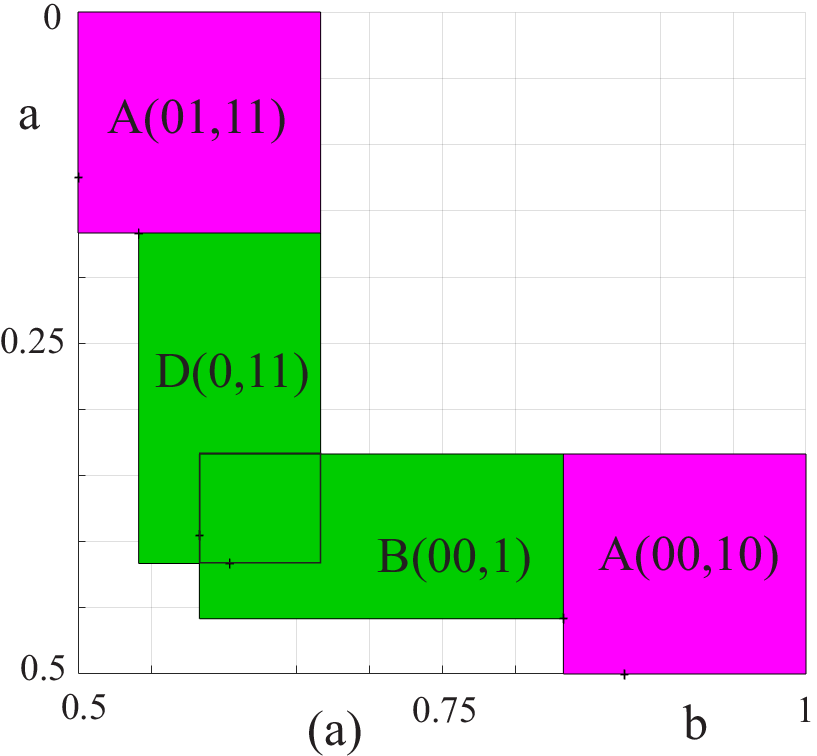}
\includegraphics[width=7cm]{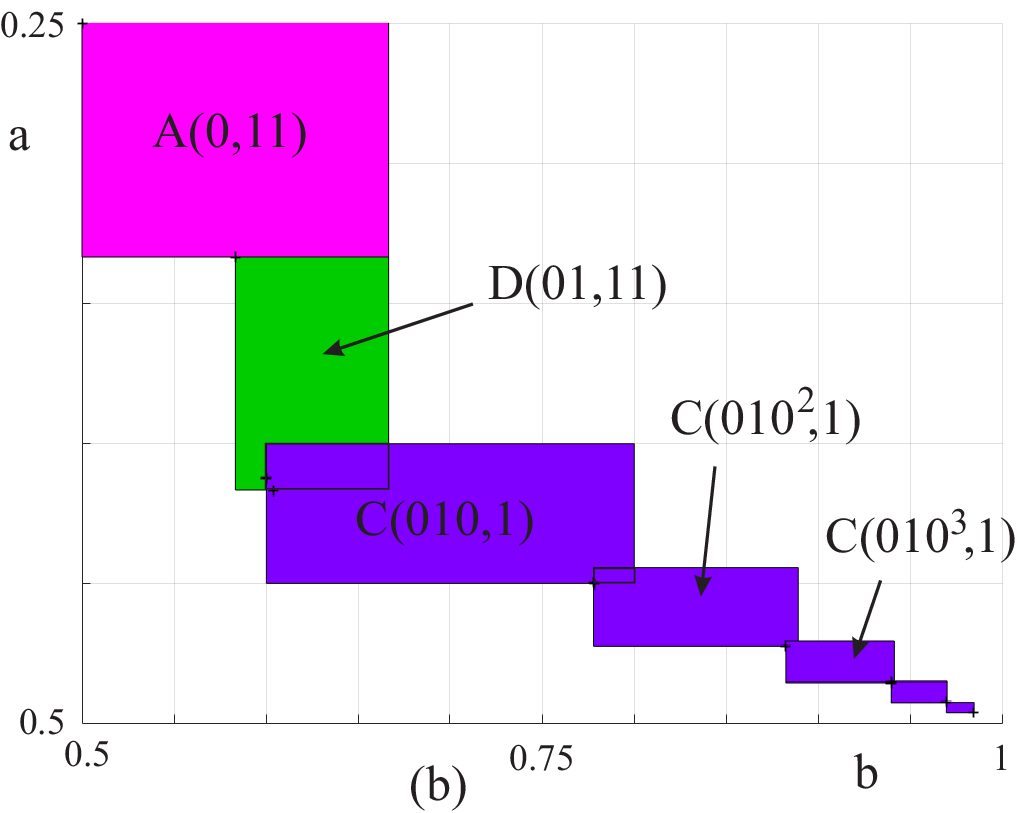}
\caption{Renormalization boxes enclosing the boundary of chaos.  (a) Class $C$, see section~\ref{sect:firstC} for details. The boxes $D(0,11)$ and $B(00,1)$ intersect on an $A$-box $A(00,11)$ which is the $A$-box of (b) and its symmetric equivalent. (b) Class $B$ see section~\ref{sect:firstB} for details.}
\label{fig:BCboxes}
\end{figure}

We can now apply these ideas to the case of complete families of plateau maps in class $C$. The decomposition of the parameter box into sub-regions with different basic induced maps is shown in Figure~\ref{fig:BCboxes}a.

\begin{theorem}\label{thm:Cind}Suppose that $F_{a,b}$ is a complete family of plateau maps in class $C$  for $(a,b)\in [a_m,a_M]\times [b_M,b_m]$. Then the boundary of chaos is contained in the finite union
\[
A(01, 11)\cup D(0,11)\cup B(00,1) \cup A(00,10).
\]
In each of the regions the induced maps form a complete family of plateau maps of the corresponding class.
\end{theorem}

\emph{Proof:} Let $a_g=\alpha (0^\infty )$, $b_g=\alpha (10^\infty)$, $a_G=\alpha (0010^\infty  )$ and  $b_G= \alpha (110^\infty )$. Then the induced map $F_1=(F^2,F)$ is a plateau map of class $B$ on $[a_g,c)\cup (c,b_g]$. Since $F(x)<\alpha (0^\infty )$ for all $x\in (\alpha (0^\infty ), \alpha (0(01)^\infty))$ the plateaus coincides with the plateaus of $F_{a,b}$ and hence $F_1$ is a complete family of plateau maps of class $B$ for $(a,b) \in B(00,1)$ where  
\[
B(00,1)=[ \alpha (0^\infty ), \alpha (0010^\infty)] \times [\alpha (110^\infty ),\alpha (10^\infty )]=[a_g,a_G]\times [b_G,b_g] .
\]
Suppose  $(a,b)\in B(00,1)$, so $F_{a,b}(x )\in [a_g,b_g]$ for all $x\in [a_g,b_g]$. If $x<a_g$ then $F(x)>a_g$ and so either $F(x)\in [a_g,b_g]$ in which case higher iterates remain in $[a_g,b_g]$ under $F_1$, or $F(x)>b_g$. Since $F(b_g)=a_g$ by definition, if $x>b_g$ then $F(b_g)<a_g$, Hence for all $x\in [a_m,b_m]$ and either there exists $n$ such that $F^n(x)\in [a_g,b_g]$ and the dynamics is determined by $F_1$, or $k(x)\in \{(01)^\infty, (10)^\infty\}$, i.e. $x$ is one of the boundary points $a_m$ or $b_m$. Thus $h_{top}(F_{a,b})>0$ if and only if $h_{top}(F_{1\, a,b})>0$. $F_1$ is a complete family of maps of class $B$ for $(a,b)\in B(00,1)$ by definition.    

A similar argument shows that $F_2=(F^2,F^2)$ on 
\[
A(01,11)=[ \alpha ([01]^\infty ), \alpha (01^\infty )] \times [\alpha (11[01]^\infty ) ),\alpha (1^\infty )]
\]
is a complete family of maps in class $A$ and $h_{top}(F_{a,b})>0$ if and only if $h_{top}(F_{2\, a,b})>0$..

The symmetric regions are 
\[
 D(0,11)= [ \alpha (01^\infty ), \alpha (001^\infty)] \times [\alpha (1101^\infty ),\alpha (1^\infty )],
\]
and
\[
A(00,10) =[ \alpha (0^\infty ), \alpha (00[10]^\infty )] \times [\alpha (10^\infty ) ),\alpha ([10]^\infty )].
\]
Using the class $C$ inequality $<_s$ of section~\ref{sect:basic},
\[
\alpha ([01]^\infty ) <  \alpha (01^\infty ) < \alpha (0^\infty ) < \alpha (001^\infty) < \alpha(0010^\infty) <\alpha (00[10]^\infty ) <\alpha (0[01]^\infty)
\]
and 
\[
\alpha ([10]^\infty ) > \alpha(101^\infty)  >  \alpha (10^\infty ) > \alpha (1^\infty ) > \alpha (110^\infty) >\alpha (11[01]^\infty ) >\alpha (1[10]^\infty)
\]
so the sets intersect as shown in Figure~\ref{fig:BCboxes}a.


It is not hard to show (the `boring' cases of \cite{G2014}) that the map has zero entropy on one side of the union of these sets (the upper right hand side as shown in Figure~\ref{fig:BCboxes}a) and positive entropy on the other. Indeed we can use the results of section~\ref{sect:bdybdy} to give a more precise description.

\begin{corollary}\label{FirstItTypeC}The boundary of chaos for a complete family of plateau maps of class $C$ is the union of continuous curves connecting the four points
\[\begin{array}{l}
P_1=(\alpha (0(01)^\infty ) ,\alpha (100(01)^\infty) ), \quad P_2=(\alpha(0010^\infty ) ,\alpha (10^\infty) ), \\ P_3=(\alpha (01^\infty),\alpha(1101^\infty ) ), \quad P_4=(\alpha(011(10)^\infty ) ,\alpha (1(10)^\infty) ).
\end{array}\]
$P_1$ is on $a=a_M$, $P_2$ is contained in $A(00,10)\cap B(00,1)$ and the curve connecting $P_1$ to $P_2$ lies entirely in $A(00,10)$; $P_3$ is contained in $D(0,11)\cap A(01,11)$ and the curve connecting $P_2$ to $P_3$ lies entirely in the region $B(00,1)\cup D(0,11)$ whose interior is a simply connected domain, and $P_4$ lies on $b=b_M$ and the curve connecting $P_3$ to $P_4$ lies in $A(01,11)$.
\end{corollary}

The corollary is a straightforward consequence of the continuous curve lemma, Lemma~\ref{lem:cc}, the lemmas of section~\ref{sect:bdybdy} and Theorem~\ref{thm:Cind}, noting that  
\begin{equation}\label{eq:listPi}
\begin{array}{rll}
P_1 & = (\alpha (0(01)^\infty ) ,\alpha (100(01)^\infty) ) & = (\alpha ([00][10]^\infty ) ,\alpha ([10][00][10]^\infty) )\\
P_2 &=(\alpha([00][10][00]^\infty ) ,\alpha ([10][00]^\infty) ) & =(\alpha([00][1][00]^\infty ) ,\alpha ([1][00]^\infty) )\\
P_3 &=(\alpha ([01][11]^\infty) , \alpha([11][01][11]^\infty ) , & =(\alpha ([0][11]^\infty), \alpha([11][0][11]^\infty ))\\
P_4 &=(\alpha(011(10)^\infty ) ,\alpha (1(10)^\infty) ) &=(\alpha ([01][11][01]^\infty) , \alpha([11][01]^\infty ).
\end{array}\end{equation}
The square brackets have been used to show the way that the points $P_i$ can be seen as the boundary points on the boundary of chaos identified in section~\ref{sect:bdybdy} for the map or its corresponding induced map. This shows that there is no discontinuity crossing the borders of the regions of different induced maps.


\section{First renormalisation for class $B$}\label{sect:firstB}
The equivalent decomposition for maps of class $B$ is a little more complicated as there are infinitely many parameter boxes covering the boundary of chaos as sketched in Figure~\ref{fig:BCboxes}b.

\begin{theorem}\label{thm:Bind}Suppose that $F_{a,b}$, $(a,b)\in [a_m,a_M]\times [b_M,b_m]$, is a complete family of plateau maps of class $B$. Then the boundary of chaos is contained in the countable union
\[
A(0, 11)\cup D(01,11)\cup \left(\bigcup_{n=0}^\infty C(010^{n+1},10^{n}) \right) .
\]
In each of the regions the indicated induced maps are basic and form a complete family of plateau maps of the corresponding class.
\end{theorem}

\emph{Proof:} The proof is essentially the same as the equivalent theorem maps of class $C$, Theorem~\ref{thm:Cind}, so in each case we simply give the region of parameter space on which the induced maps form a complete family. The infinite intersection has simply connected interior so there is no issue about proving continuity across the boundaries of those parameter regions, so the only consistency conditions that need to be established is that of showing that the boundary of chaos is continuous at the intersection of the boundaries of $A(0,11)$ and $D(01,11)$ and on $b=b_M$.

By Lemma~\ref{lem:posbdy} and the symbolic remarks immediately after its proof for case $D$, the intersection of the boundary of chaos for $A(0,11)$ are
\begin{equation}\label{eq:bdyptsA}
\begin{array}{rll}
P_1 & =(\alpha ([0][11][0]^\infty ,\alpha ([11][0]^\infty ) & \\ 
P_2 & =(\alpha ([0][11]^\infty ), \alpha ([11][0][11]^\infty )) &= (\alpha ([01][11]^\infty ), \alpha ([11][01][11]^\infty )).
\end{array}\end{equation}
The first of these is the limit from above of the boundary of entropy for type-$B$ maps and the second, in the second re-arrangement, is the boundary of entropy for $D(01,11)$ maps. The limit of the $C$-boxes is clearly $(010^\infty ,10^\infty)$, the other end point for the boundary of chaos on the boundary of parameters for maps of class $B$ by the comments after Lemma~\ref{lem:posbdy}.
 

\section{Boundary of chaos for class $A$}\label{sect:bdytypeA}
This case has been studied by many people. The key results are from Gambaudo et al \cite{GPTT} for maps (though without any discussion of full families) and Glendinning and Sidorov \cite{GS2015} for open maps. In this section we rephrase these results in the same language as the previous sections and go on to characterise the boundary of chaos. 

Recall that $r^\pm_{p/q}$ are balanced words (Definition~\ref{def:balanced}).

\begin{figure}
\centering
	\includegraphics[width=9cm]{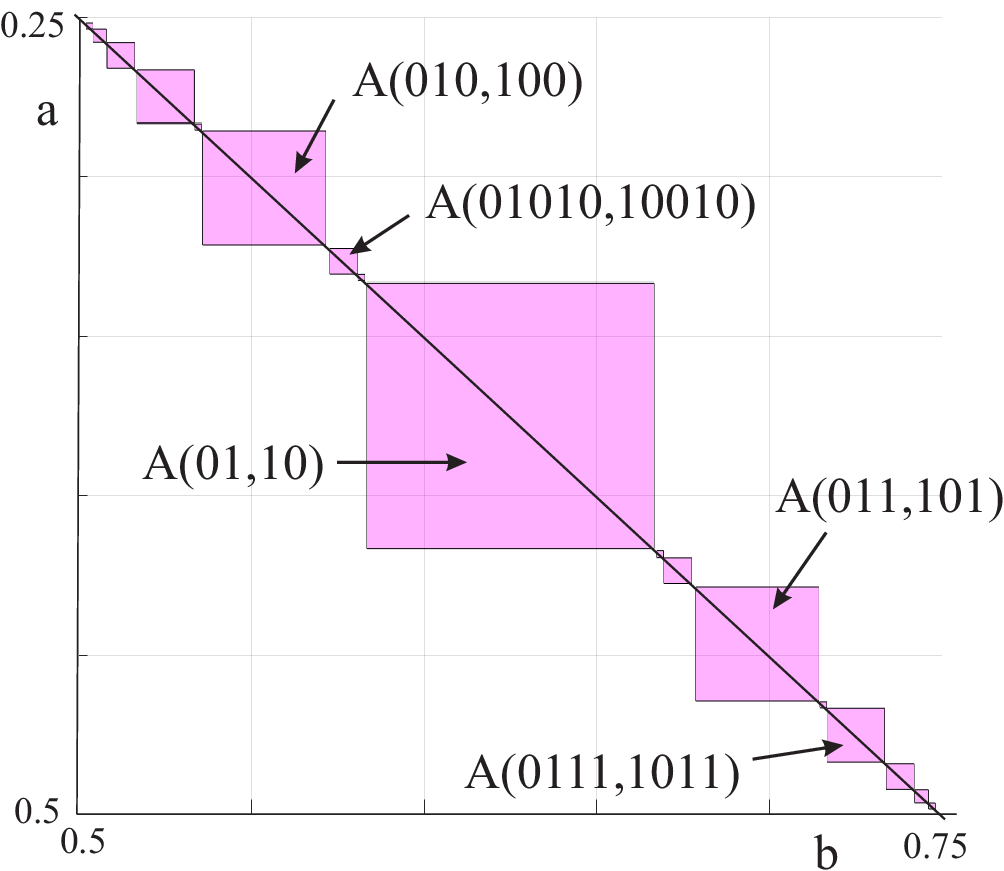}
	\caption{Renormalization boxes for class $A$ and the line $b-a=0.25$.}
	\label{fig:Aboxes}
\end{figure}

For each $p/q\in \mathbb{Q}$ there is a region $A(r^-_{p/q},r^+_{p/q})$ for which the corresponding induced map $G=(F^q,F^q)$ is well-defined and $F$ has positive entropy if and only if $G$ has positive entropy, e.g.  \cite{G1990}. An important feature of type $A$ maps is the existence of a curve $R$ (Figure~\ref{fig:Aboxes}.) in parameter space on which $F^2(a)=F^2(b)$ with $F(b)\le c\le F(a)$. If $(a,b)\in R$ then  the map restricted to $[F(b),F(a)]$ can be seen as a continuous monotonic circle map and so all the standard results for circle maps apply \cite{GPTT}. In particular \cite{GS2015,Veerman1989} there is set of Hausdorff dimension zero on which the dynamics has irrational rotation number and there are no periodic orbits, and all other parameter values lie in the sets $A(r^-_{p/q},r^+_{p/q})$. 

    
Theorem~\ref{FirstItTypeA} and Corollary~\ref{cor:++} below are the results of Glendinning and Sidorov \cite{GS2015} re-interpreted in the language of induced maps. More specifically, Corollary~\ref{cor:++} corresponds to Theorem~2.13 of \cite{GS2015}. We state these results without further proof.

\begin{theorem}\label{FirstItTypeA}Suppose that $F_{a,b}$, $(a,b)\in [a_m,a_M]\times [b_M,b_m]$, is a complete family of plateau maps of class $A$ derived from an expanding double cover. Then the boundary of chaos is contained in the union
\[
\Omega = \left(\bigcup_{{p<q,~ p,q ~coprime}}A(r^-_{p/q},r^+_{p/q})\right) \cup I_\omega
\]
where $I_\omega$ is the union of sets on which the dynamics is conjugate to a circle map with irrational rotation number. On each of the regions $A(r^-_{p/q},r^+_{p/q})$ the corresponding induced map is a complete family of plateau maps of class $A$.
\end{theorem}

Let $A(r^-_{p/q},r^+_{p/q})=A_{p/q}$ and let 
\begin{equation}
{\mathcal S}=R\backslash \left(\bigcup_{{p<q,~ p,q ~coprime}}{\rm int}A(r^-_{p/q},r^+_{p/q})\right).
\label{eq:cals}
\end{equation}

\begin{corollary}Suppose that $F_{a,b}$, $(a,b)\in [a_m,a_M]\times [b_M,b_m]$, is a complete family of plateau maps of class $A$. Then the boundary of chaos is a continuous curve $\mathcal{C}$ and if $(a,b)\in \mathcal{C}$ and either
\begin{itemize}
\item[(i)]  there exists a finite (possibly empty) sequence $p_i/q_i$, $i=1, \dots ,n$ such that $(a,b)$ is in the $A_{p_n/q_n}$ box of the $A_{p_{n-1}/q_{n-1}}$ box of the $\dots$  $A_{p_1/q_1}$ box, and the boundary points are in the equivalent of the set ${\cal S}$ of (\ref{eq:cals});
\item[(ii)] there exists a finite sequence $p_i/q_i$, $i=1, \dots ,n$ such that $(a,b)$ is in the $A_{p_n/q_n}$ box of the $A_{p_{n-1}/q_{n-1}}$ box of the $\dots$  $A_{p_1/q_1}$ box, and either $(k_-,k_+)=(v^\infty,uv^\infty)$, $(k_-,k_+)=(vuv^\infty,uv^\infty)$, $(k_-,k_+)=(vu^\infty,u^\infty)$ or $(k_-,k_+)=(vu^\infty,uvu^\infty)$;
\item[(iii)] there exists an infinite sequence $p_i/q_i$, $i=1, 2,\dots $ such that for all $n> 1$ $(a,b)$ is in the $A_{p_n/q_n}$ box of the $A_{p_{n-1}/q_{n-1}}$ box. 
\end{itemize}
Case (ii) occurs on line segments in parameter space, cases (i) and (iii) contain no connected arcs. 
\label{cor:++}
\end{corollary}

Note that here case (ii) is given in terms of the kneading sequences whilst in \cite{GS2015} it is described through the addresses of $a$ and $b$. The only non-standard part of this result is the final comment. Irrational rotations occur at points on the curve \cite{Veerman1989} so for each irrational number (i) is a single point. Case (ii) occurs on lines segments by Theorem~\ref{thm:hettrans} or \cite[Theorem 2.13]{GS2015} which include the `boundary on the boundary' points for case A of section~\ref{sect:bdybdy} as end points, and case (iii) on points using the expanding property $\lambda >1$ in Definition~\ref{def:plateau} which implies that the size of successive parameter boxes decreases exponentially (this argument is expanded on in more detail in the proof of Theorem~\ref{thm:main} of the next section).

Smooth non-invertible circle maps and and non-expanding piecewise smooth maps such as (\ref{eq:smoothcases}) with $a, b >0$  have a codimension one transition which does not occur in the plateau maps defined here. In the cases studied here the existence of the well-ordered $p/q$ periodic orbit is defined by $a\ge \alpha((r^-_{p/q})^\infty)$ and $b\le\alpha((r^+_{p/q})^\infty)$ in the notation of section~\ref{sect:chaos}. Thus the periodic orbits exist in a region defined by vertical and horizontal lines.  Consider the period two renormalization $A_{1/2}$. The boundaries not defined $\alpha ((r^\pm_{1/2})^\infty)$ ensure that the induced maps $(F^2,F^2)$ map an interval into itself for parameters in $A_{1/2}$ and so orbits in this interval have even periods. One way to understand the lower boundaries of $A_{1/2}$ (i.e. with $b<a+0.25$) is that $A$-boxes with $\frac{p}{q}\to \frac{1}{2}$ tend to $A_{1/2}$ and intersect $b-a=0.25$ on intervals. Since their upper boundaries for the existence of the well-ordered orbits are vertical and horizontal, their accumulation is also on vertical and horizontal lines, which must accumulate on the vertical and horizontal boundaries of $A_{1/2}$ through the intersection of $A_{1/2}$ with $b-a=0.25$. By contrast, in the smooth cases the accumulation of these saddle node curves does not extend into the equivalent of $A_{p/q}$ immediately after the loss of invertibility, and it is this interval on the saddle node curve that provides the codimension one circle intermittency transition to chaos in those cases.      


\section{Boundary of chaos for classes $B$, $C$ and $D$}\label{sect:mainresult}
The same argument following the possible basic renormalizations can be used to characterise the boundary of chaos for the other classes of maps.

\begin{theorem}\label{thm:main}Suppose that $F_{a,b}$ $(a,b)\in [a_m,a_M]\times [b_M,b_m]$ is a complete family of plateau maps of class $B$, $C$ or $D$ derived from an expanding double cover. Then the boundary of chaos is a continuous curve $\mathcal{C}$ and if $(a,b)\in \mathcal{C}$ either
\begin{itemize}
\item[(i)] after a finite sequence of $B$, $C$ and $D$ renormalizations $(a,b)$ is in an $A$ box and one of the three possibilities of Corollary~\ref{cor:++} holds;
\item[(ii)] after a finite sequence of $B$, $C$ and $D$ renormalizations there is an infinite sequence of $B(00,1)$ then $C(010,1)$ renormalizations or an infinite sequence of $D(0,11)$ then $C(0,101)$ renormalizations; 
\item[(iii)] there is an infinite sequence of $B$, $C$ and $D$ renormalizations not in case (ii).    
\end{itemize}
Case (ii) occurs on line segments in parameter space and (iii) at points.
\end{theorem}

\emph{Proof:} 
In Corollary \ref{FirstItTypeC}, Theorem \ref{thm:Bind} and Theorem \ref{FirstItTypeA}, we defined sets of boxes that contain the boundary of chaos and points $P_i$ that belong to the boundary of chaos. In each box of these sets, we define similar sets of boxes and prove with the same method that these sets all contain a part of the boundary of chaos. These different boxes within a box are connected to each other through the points $P_i$. We then use an induction to prove that after an infinite number of boxes within boxes, the final set contains the entire boundary of chaos.

The only part left that requires proof is now the statement about whether limits are points or line segments. This will prove that the final set defined above is contained in the boundary of chaos and therefore is equal to the boundary of chaos.

Each renormalization at level $n$ reduces the size of the $a$-interval of the box-within-a-box structure by a factor $\lambda_n^{-(|\omega_-|}$, where $\lambda_n$ is the smallest slope of the induced map at the previous level. Similarly for the $b$-interval. Thus if there are $k$ renormalizations with $\min (|\omega_-|, |\omega_+|)\ge 2$ the length of the sides of the induced boxes in parameter space are both less than $(b_m-a_m)\lambda^{-2^k}$.  If there are an infinite number of these then the size of each side of the boxes tends to zero and hence the limit is a point.

The renormalizations of classes $B$, $C$ and $D$ which do not have $\min (|\omega_-|, |\omega_+|)\ge 2$ are $D(0,11)$ and $B(00,1)$ in a type $C$ renormalization (Theorem~\ref{thm:Cind}), and $C(010,1)$ in type B (Theorem~\ref{thm:Bind}) and $C(0,101)$  in type D. If $C(010,1)$ is followed by a $D(0,11)$ then the net effect is to increase the number of iterates on both sides of the map and a further factor of at least $\lambda_n^{-1}$ reduces the length of each side, similarly if $C(0,101)$ is followed by $B(00,1)$. Therefore to avoid both sides going to zero the sequence must eventually be $C(0,101)$ followed by $D(0,11)$ followed by $C(0,101)$ and so on, or similarly with $C(010,1)$ and $B(00,1)$. In this case the slope of one branch of the map goes to infinity (so one side of the iterated boxes goes to zero) whilst the other is unchanged, and hence in the limit the boxes on which these induced maps are defined tends to a line segment. 
\newline\rightline{$\square$}


\begin{figure}
\centering
	\includegraphics[width=4cm]{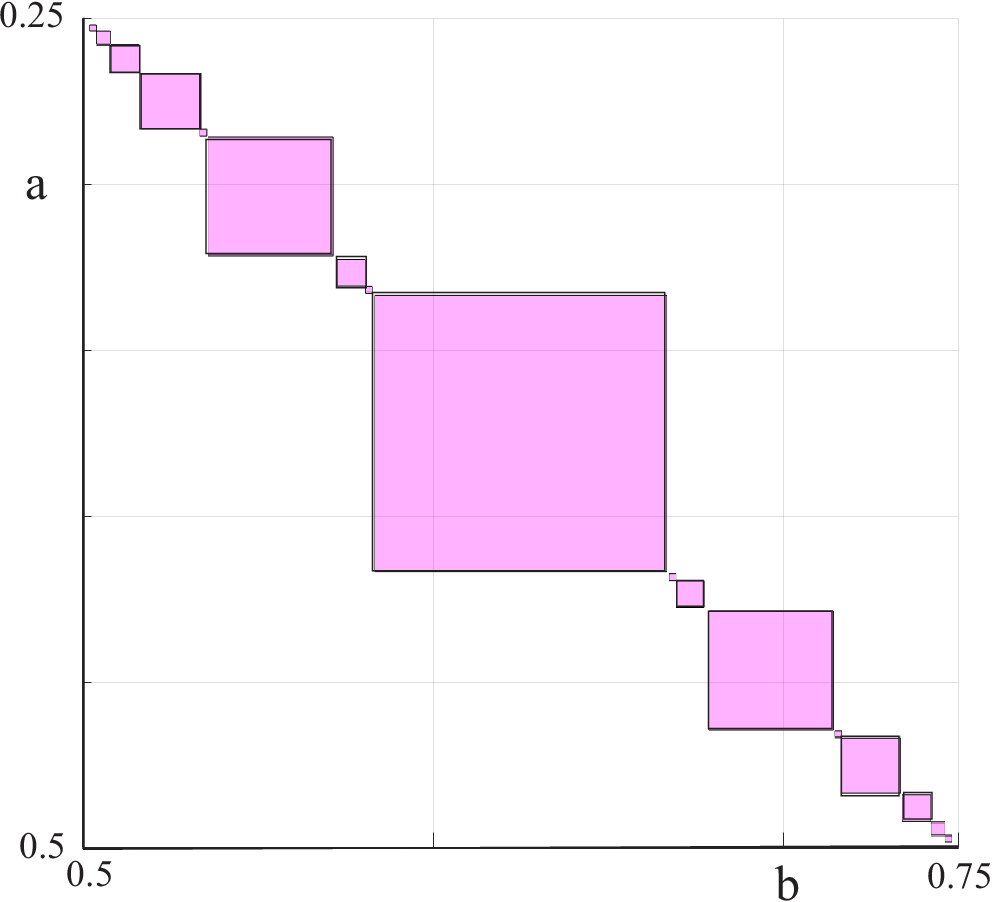}
\hspace{0.5cm}
\includegraphics[width=4cm]{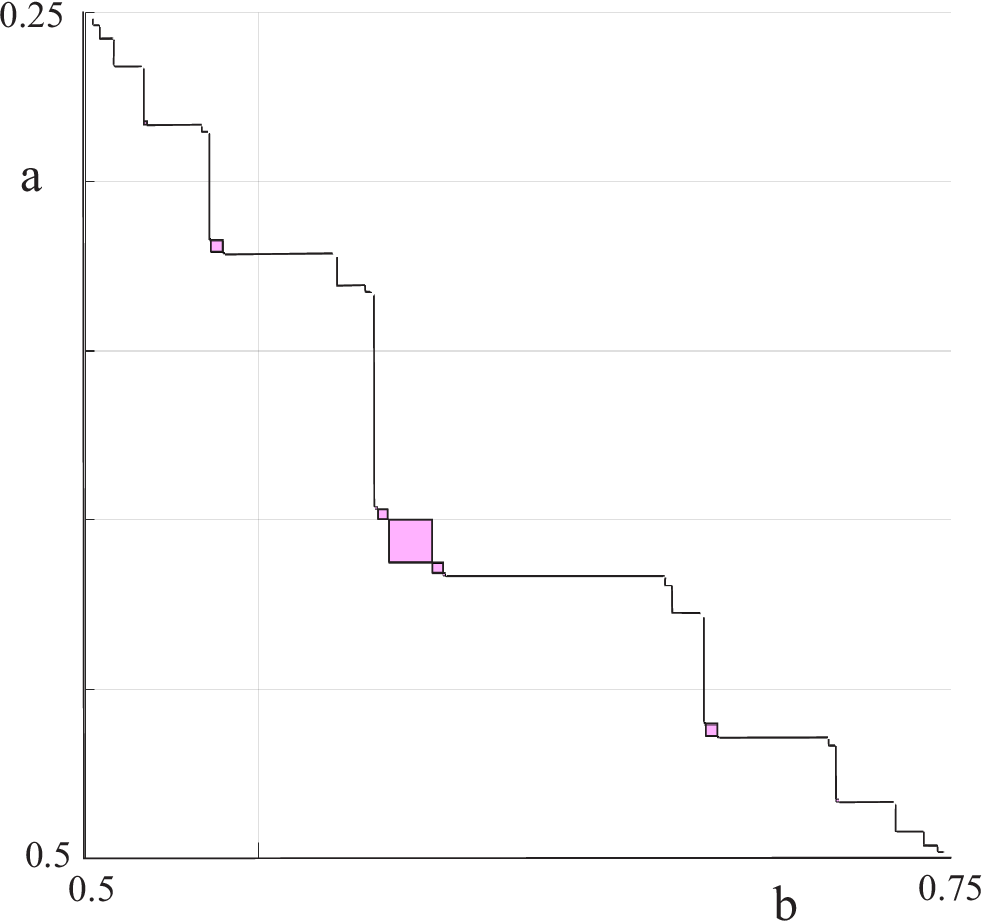}
\hspace{0.5cm}
\includegraphics[width=4cm]{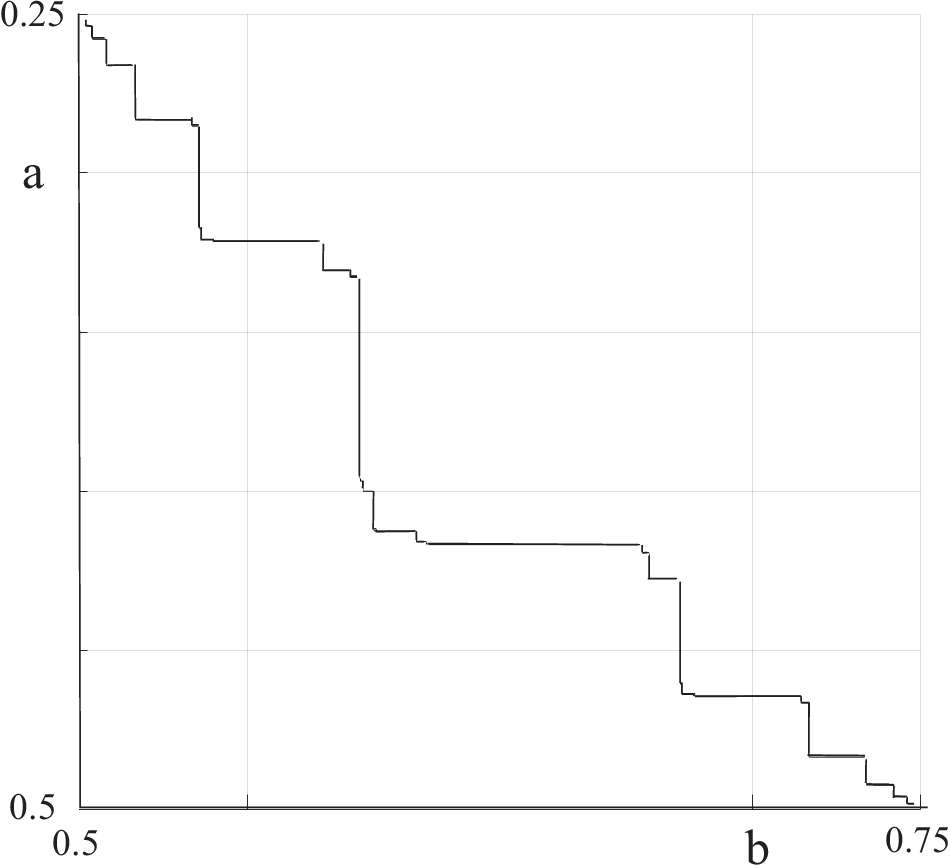}\\
\vspace{0.3cm}
	\includegraphics[width=4cm]{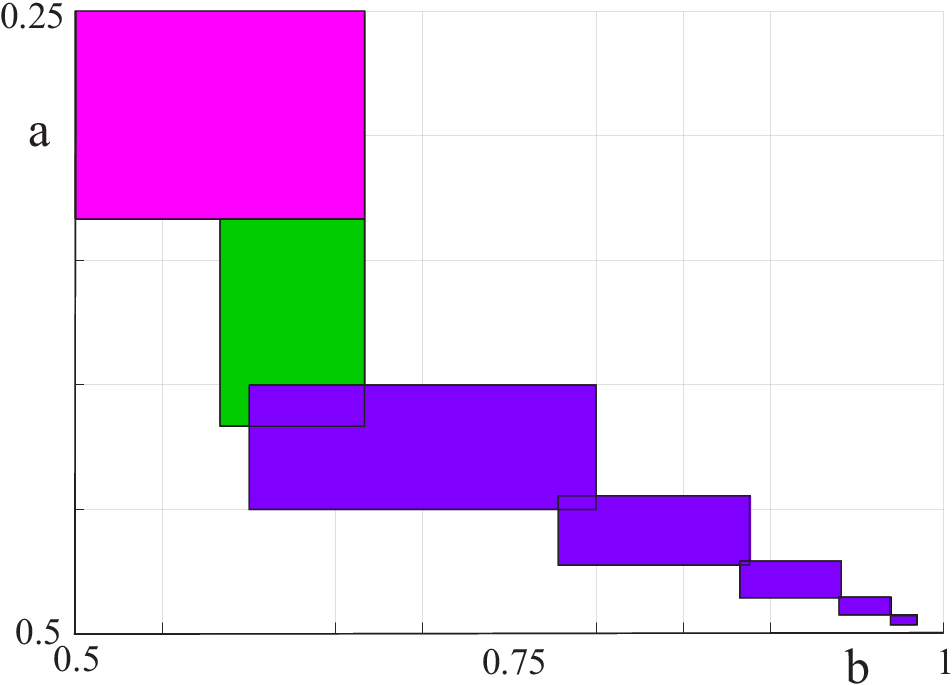}
\hspace{0.5cm}
\includegraphics[width=4cm]{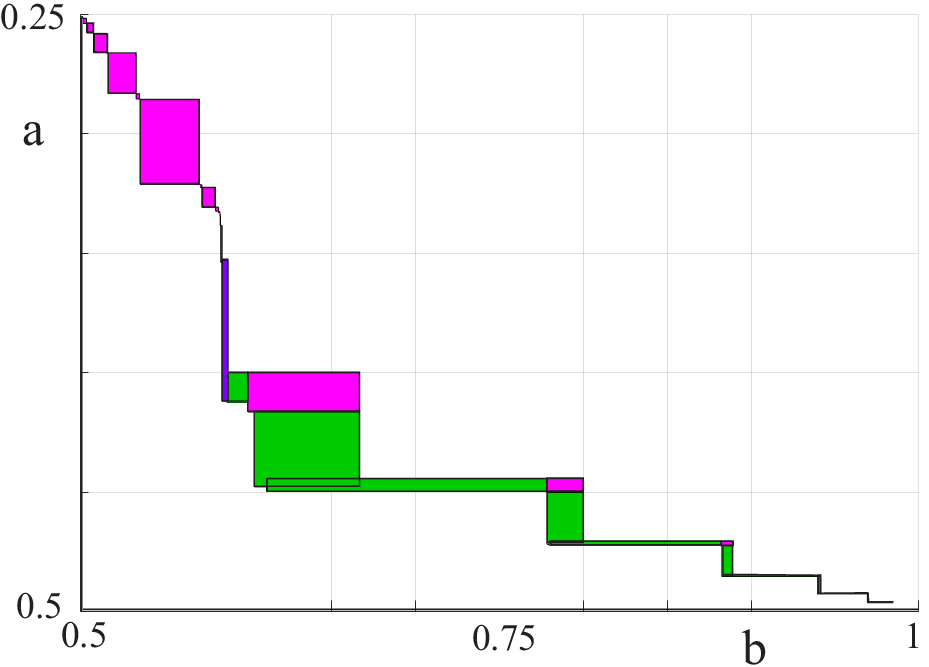}
\hspace{0.5cm}
\includegraphics[width=4cm]{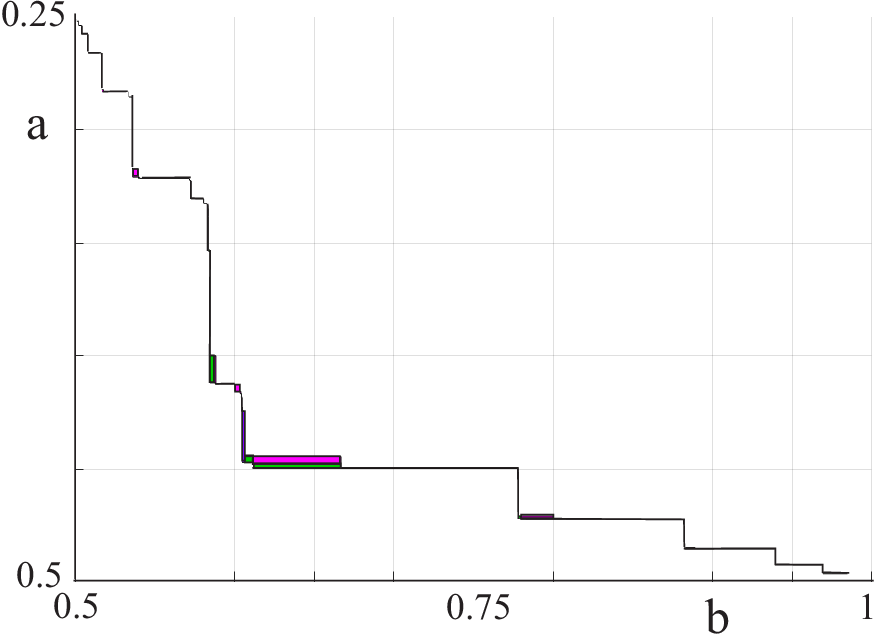}\\
\vspace{0.3cm}
	\includegraphics[width=4cm]{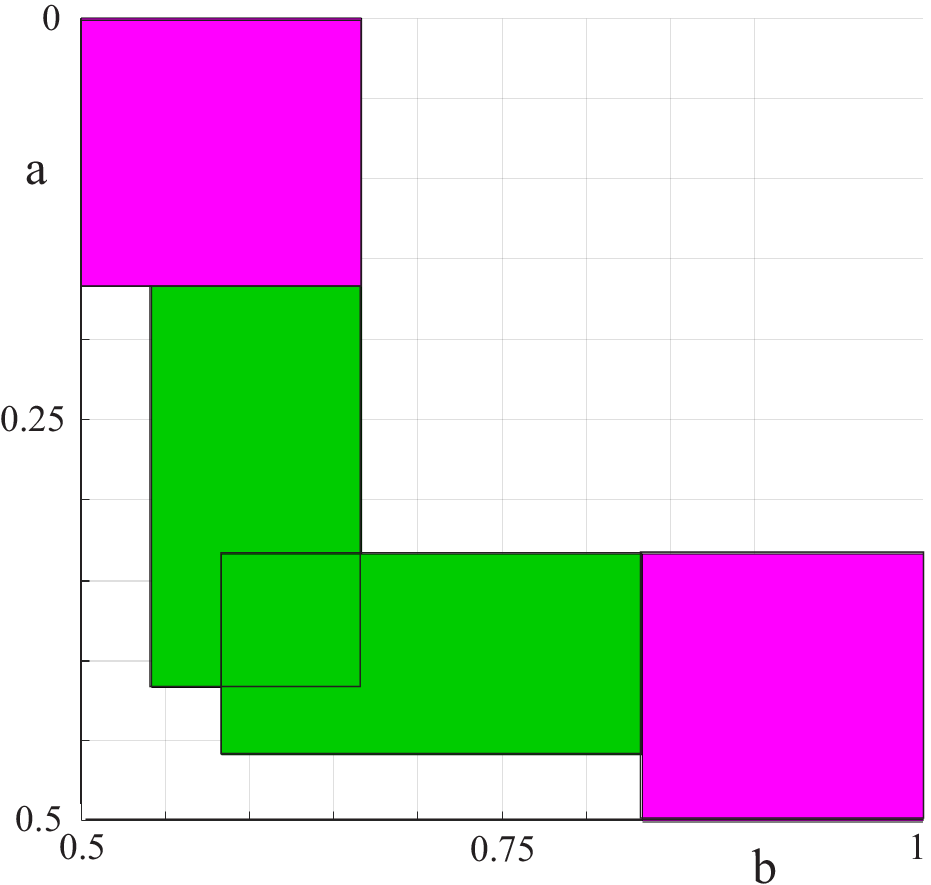}
\hspace{0.5cm}
\includegraphics[width=4cm]{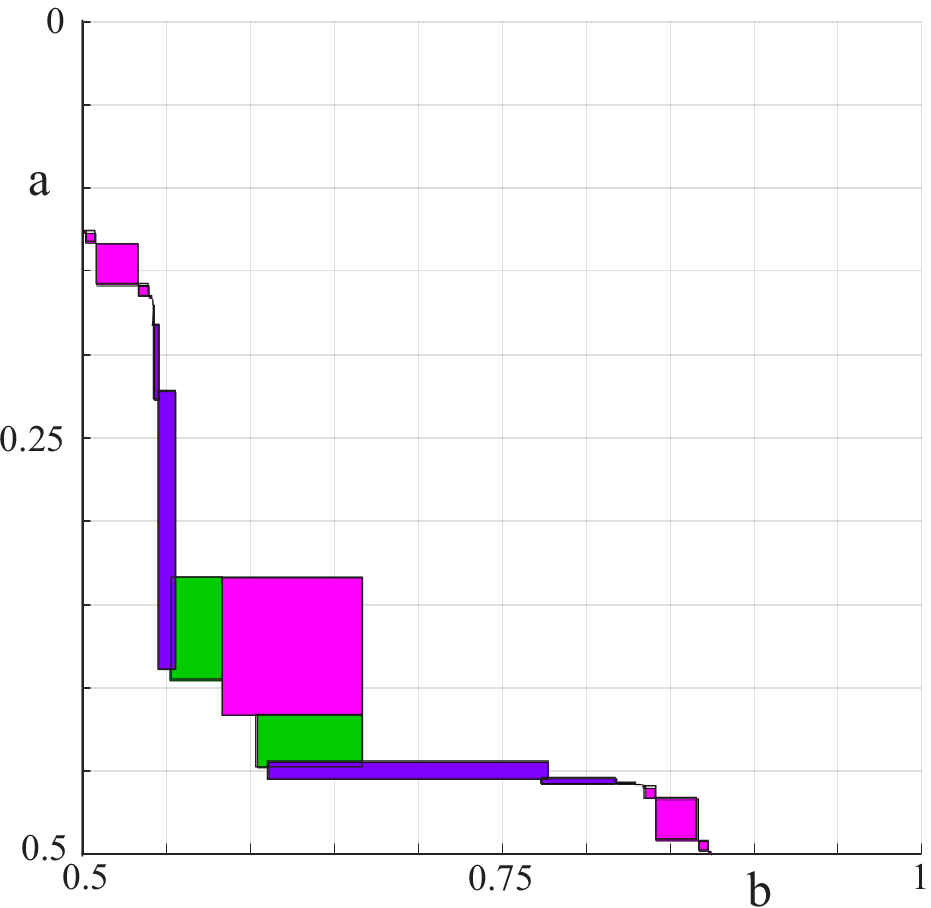}
\hspace{0.5cm}
\includegraphics[width=4cm]{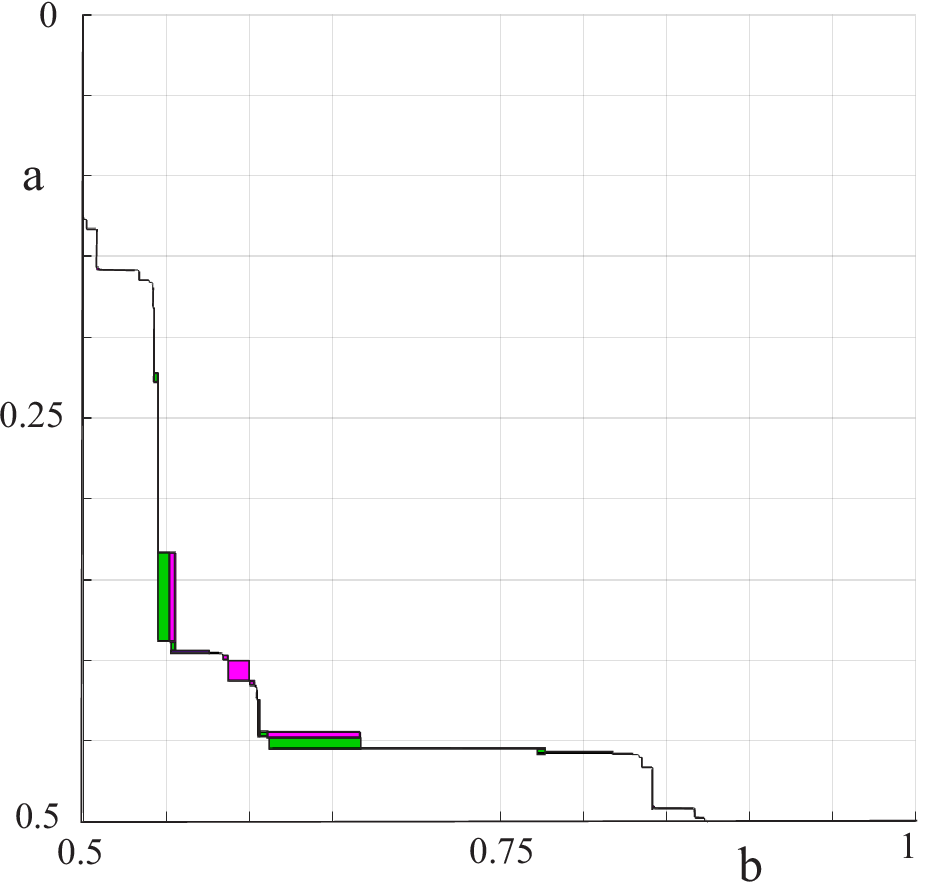}
\caption{Three iterations of the inductive process to construct the boundary of chaos. First row, class A; second row, class B; and third row class C. The first diagrams in each row are the parameter regions of Figure~\ref{fig:Aboxes} and Figure~\ref{fig:BCboxes} respectively.} 
	\label{IterationsBoxesGraph}
\end{figure}


This iterative process is illustrated in Figure~\ref{IterationsBoxesGraph} for the case when the modulus of the slopes of the original double covers equals two.  Note that case (ii) of Theorem~\ref{thm:main} is the anharmonic route described by Glendinning \cite{G1993} and case (iii) includes (asymptotically) the standard Morse-Thue renormalization, with $B$ and $D$ renormalizations alternating \cite{AS1999}.


\section{Discussion}\label{sect:conclusion}We have given a complete description of the boundary of chaos for plateau maps derived from expanding double covers. Since the complete families considered are also full, in the sense that every consistent kneading invariants are realised by an element of the family, this describes all the possible transitions for piecewise smooth maps with two monotonic branches. An important conclusion is that there are only two codimension one transitions: the heteroclinic transition of \cite{MT1986} and the anharmonic cascade of \cite{G1992,G1993}. 

The correspondence between the property of positive entropy for plateau maps and positive Hausdorff dimension of survivor sets of some open maps means that this discussion completes a programme to characterise the boundary between positive and zero Hausdorff dimension in these maps. It also shows that whilst the codimension one components of the boundary in the monotonic increasing case considered by \cite{GS2015} have only a finite number of periodic orbits, there are codimension one segments in the more general cases for which there is an infinite number of periodic orbits on the boundary.

The anharmonic cascade replaces the more well-known period-doubling cascade as the robust transition with infinite bifurcations. Whilst the period-doubling cascade leads to the Morse-Thue sequence \cite{AS1999}, it is not clear how the symbol sequences for that generate the anharmonic cascade, (\ref{eq:replaceC}) and (\ref{eq:replaceB}), relate to known sequences. Given their central role in this natural context it would be interesting to determine whether they feature in other studies.


\end{document}